\documentclass[11pt,twoside]{amsart}
\usepackage{amssymb,amsmath}
\usepackage{amscd}
\usepackage{amsthm}
\usepackage{latexsym}
\usepackage[noadjust]{cite}
\usepackage{mathrsfs}
\usepackage{indentfirst}
\usepackage{graphicx}
\usepackage{graphicx}
\usepackage{subfigure}
\usepackage{enumitem}
\allowdisplaybreaks[4]

\def\dfrac{\displaystyle\frac}
\let\oldsection\section
\renewcommand\section{\setcounter{equation}{0}\oldsection}

\newtheorem{theorem}{Theorem}[section]
\newtheorem{lemma}[theorem]{Lemma}

\newtheorem{remark}[theorem]{Remark}
\newtheorem{example}{Example}
\newtheorem{corollary}[theorem]{Corollary}

\usepackage[hidelinks,
hypertexnames=false,
linktocpage,pdfpagelabels,
bookmarksnumbered,bookmarksopen]{hyperref}
\usepackage[hyperpageref]{backref}

\makeatletter
\newcommand{\itemcustom}[2]{%
	\refstepcounter{enumi}%
	\item[#1]%
	\protected@edef\@currentlabel{#1}%
	\label{#2}%
}
\makeatother

\makeatletter
\providecommand\@dotsep{5}
\def\listtodoname{List of Todos}
\def\listoftodos{\@starttoc{tdo}\listtodoname}
\makeatother

\newcommand{\vp}{\varphi}

\newcommand{\eps}{\varepsilon}

\def\N{\mathbb{N}}
\def\R{\mathbb{R}}

\def\D{\mathcal{D}}

\def\RN{\mathbb{R}^N}

\def\n{\nabla}

\def\a{\alpha}
\def\b{\beta}
\def\l{\lambda}
\def\g{\gamma}

\def\t{\theta}
\def\irn{\int_{\RN}}

\newcommand{\weakto}{\rightharpoonup}

\newcommand{\cC}{{\mathcal C}}

\newcommand{\cE}{{\mathcal E}}

\newcommand{\cM}{{\mathcal M}}

\newcommand{\cQ}{{\mathcal Q}}

\newcommand{\cS}{{\mathcal S}}

\title[Normalized Born--Infeld ground states]{Normalized ground state solutions to the nonlinear Born--Infeld problem}

\author[J. Mederski]{Jaros\l aw Mederski}
\author[X. Zeng]{
	Xiangjian Zeng}

\address{J. Mederski
	\newline\indent Institute of Mathematics,
	\newline\indent Polish Academy of Sciences  
	\newline\indent ul. \'Sniadeckich 8, 00-656 Warsaw, Poland
	\newline\indent and
	\newline\indent Faculty of Mathematics and Computer Science,
	\newline\indent Nicolaus Copernicus University,
	\newline\indent ul. Chopina 12/18, 87-100 Toru\'n, Poland
}
\email{\href{mailto:jmederski@impan.pl}{jmederski@impan.pl}}

\address{X. Zeng
	\newline\indent Institute of Mathematics,
	\newline\indent Polish Academy of Sciences  
	\newline\indent ul. \'Sniadeckich 8, 00-656 Warsaw, Poland
	\newline\indent and
	\newline\indent Faculty of Mathematics, Informatics and Mechanics 
	\newline\indent University of Warsaw 
	\newline\indent ul. Banacha 2, 02-097 Warsaw, Poland}
\email{\href{mailto:xzeng@impan.pl}{xzeng@impan.pl}}


\begin{document}
	\begin{abstract}
		We study normalized solutions to the Born--Infeld problem
		\begin{equation*}
			-{\rm div}\left(\frac{\n u}{\sqrt{1-|\n u|^2}}\right)+\mu u=f(u)\ \hbox{ in }\RN,\qquad \irn|u|^2\,dx=2m^2,
		\end{equation*}
		$N\ge3$, where $m>0$ is prescribed and $\mu>0$ is an unknown Lagrange multiplier. The nonlinearity $f$ is assumed to be superlinear at infinity, with no upper restriction on its growth, and to behave like $|s|^{p-2}s$ near zero, with $2<p<2^*$. Following a Lagrangian approach, we relate the action ground state level $a(\mu)$ of the prescribed frequency problem to the masses of its ground states by means of one-sided derivative formulas, which lead to the existence of normalized solutions. A distinctive feature of the Born--Infeld operator is that the gradient constraint makes the energy bounded from below on every mass sphere, even in the $L^2$-supercritical case. For $2+4/N\le p<2^*$ we obtain a threshold $c_*>0$ such that normalized action ground states exist for $m^2>c_*$, and none exist for $m^2<c_*$. In the $L^2$-subcritical case they exist for every $m^2>0$. We prove that every action ground state has constant sign and is radial up to a translation.
	\end{abstract}
	\maketitle
	
	\section{Introduction}\label{sec intro}
	
	The Born--Infeld theory \cite{Bnat,B,BInat,BI} was proposed as a nonlinear alternative to Maxwell's electrodynamics in which the field energy of a point charge is finite. In the electrostatic case the potential solves an equation involving the operator $-{\rm div}(\n u/\sqrt{1-|\n u|^2})$, and the constraint $|\n u|<1$ expresses the fact that the electric field is bounded, the absolute field constant of Born being normalized to one. The same operator gives the mean curvature of spacelike graphs in the Lorentz--Minkowski space \cite{BaS}. Even in the electrostatic case with a prescribed charge density, the regularity of the minimizers of the Born--Infeld energy, and in particular the question whether the gradient bound is strict, is a delicate issue, see Bonheure and Iacopetti \cite{Bon-Iac2,Bon-Iac3}, and Byeon et al{.}\ \cite{BIMM2024}. The variational density of the operator
	\[
	j_{\rm BI}(r)=1-\sqrt{1-r^2}
	=\frac{r^2}{2}+\frac{r^4}{8}+O(r^6)\qquad(r\to0)
	\]
	behaves like the Dirichlet density for small gradients, but it is not homogeneous and its derivative blows up as $r\to1^-$. These two features, the quadratic behavior for small gradients and the bound on the admissible gradients, are responsible for the phenomena studied in this paper, which is devoted to the normalized problem
	\begin{equation}\label{eq}
		\begin{cases}
			-{\rm div}\left(\dfrac{\n u}{\sqrt{1-|\n u|^2}}\right) +\mu u = f(u) \qquad \text{in } \R^N, \quad N\geq 3,\\[5mm]
			\irn |u|^2\, dx =2 m^2,
		\end{cases}
	\end{equation}
	where $m>0$ is prescribed and $\mu>0$ is an unknown Lagrange multiplier. Two questions are addressed: which masses are attained by the action ground states of the prescribed frequency problem, and which of these ground states minimize the energy on the mass sphere. The two questions have different answers: in the $L^2$-supercritical case we find two different mass thresholds, one for each notion of ground state.
	
	The time-dependent electromagnetic Born--Infeld equations have the structure of a system of conservation laws. Brenier \cite{Brenier} showed that, once the energy density and the Poynting vector are added to the unknowns, they form a hyperbolic system of ten conservation laws with a convex entropy which resembles the equations of magnetohydrodynamics. In the low field regime this augmented system reduces to the Maxwell equations and in the high field regime to pressureless magnetohydrodynamics. For variational results in the presence of a magnetic field we refer to \cite{dAM}. The scalar problem \eqref{eq} admits a different, formal Hamiltonian interpretation, which we now describe.
	
	In nonlinear Schr\"odinger models the prescribed $L^2$-norm represents the number of particles or the power of a beam, see \cite{Caz,LSY2000}. In our scalar setting, a formal mass-preserving evolution is obtained by extending the energy to complex-valued functions $\psi=q+ip$. Setting
	\[
	\mathscr H(\psi)=\irn\left(1-\sqrt{1-|\nabla\psi|^2}-F(|\psi|)\right)\,dx,
	\qquad |\nabla\psi|^2=|\nabla q|^2+|\nabla p|^2,
	\]
	where $F(s)=\int_0^sf(t)\,dt$, the canonical Hamilton equations $q_t=\delta\mathscr H/\delta p$, $p_t=-\delta\mathscr H/\delta q$ take the form
	\[
	i\partial_t\psi
	=-{\rm div}\left(\frac{\nabla\psi}{\sqrt{1-|\nabla\psi|^2}}\right)
	-\frac{f(|\psi|)}{|\psi|}\psi,
	\]
	where the nonlinear term is extended by continuity at $\psi=0$. If $f$ is odd, then $\mathscr H$ coincides on real-valued functions with the energy $\cE$ introduced in Subsection \ref{subsec framework}. The Hamiltonian is invariant under the phase shifts $\psi\mapsto e^{i\vartheta}\psi$, and for smooth, sufficiently decaying solutions with $\|\nabla\psi\|_\infty<1$ a direct computation gives the local conservation law
	\[
	\partial_t\frac{|\psi|^2}{2}
	+\operatorname{div}\left(
	\frac{\operatorname{Im}(\overline\psi\,\nabla\psi)}
	{\sqrt{1-|\nabla\psi|^2}}\right)=0,
	\]
	hence $\frac12\|\psi\|_2^2$ is conserved. If $f$ is odd, the ansatz $\psi(t,x)=e^{i\mu t}u(x)$ leads to the equation in \eqref{eq}, and the constraint in \eqref{eq} prescribes the conserved quantity. The same conclusion holds for the second order evolution
	\begin{equation}\label{BI wave equation}
		\partial_{tt}\phi-{\rm div}\left(\frac{\nabla\phi}{\sqrt{1-|\nabla\phi|^2}}\right)+\kappa^2\phi=\frac{f(|\phi|)}{|\phi|}\phi,\qquad\kappa>0,
	\end{equation}
	which replaces the Laplacian in the nonlinear Klein--Gordon equation by the Born--Infeld operator and is the Hamiltonian system associated with the energy
	\begin{equation*}
		\irn\Big(\frac12|\partial_t\phi|^2+1-\sqrt{1-|\nabla\phi|^2}+\frac{\kappa^2}{2}|\phi|^2-F(|\phi|)\Big)\,dx
	\end{equation*}
	of a complex-valued field $\phi$. The invariance under the phase shifts gives the conserved charge ${\rm Im}\irn\overline\phi\,\partial_t\phi\,dx$, and, if $f$ is odd, the standing waves $\phi(t,x)=e^{i\omega t}u(x)$ with $0<\omega<\kappa$, the analogues of the $Q$-balls of \cite{Col,LP}, solve \eqref{eq-mu} with $\mu=\kappa^2-\omega^2$ and have charge $\omega\lVert u\rVert_2^2=2\omega m^2$. Thus prescribing the charge at a given frequency amounts to prescribing the mass in \eqref{eq}. We stress that these are only formal motivations for the mass constraint, and not derivations from the electromagnetic system studied in \cite{Brenier}. 
	
	A further motivation comes from the geometry of spacelike hypersurfaces
	in the Lorentz--Minkowski space $\mathbb L^{N+1}$, that is, $\R^{N+1}$ endowed with the
	metric $dx_1^2+\cdots+dx_N^2-dz^2$. The connection between the prescribed
	mean curvature equation in $\mathbb L^{N+1}$ and the Born--Infeld variational problem is developed
	in \cite{BaS,BdAP,Bon-Iac2,Bon-Iac3}. For $u\in\cC^2(\RN)$ with
	$|\nabla u|<1$, the graph $\Gamma_u=\{(x,u(x)):x\in\RN\}$ is
	spacelike and its induced area element is
	$\sqrt{1-|\nabla u|^2}\,dx$, and the Born--Infeld density $1-\sqrt{1-|\nabla u|^2}$ therefore measures
	the local area deficit of $\Gamma_u$ with respect to a horizontal hyperplane.
	If the mean curvature of $\Gamma_u$ is normalized as
	\[
	H_u=\frac1N\operatorname{div}
	\left(\frac{\nabla u}{\sqrt{1-|\nabla u|^2}}\right),
	\]
	then \eqref{eq} is a prescribed mean curvature problem in which the
	curvature is a function of the height, namely $H_u=(\mu u-f(u))/N$,
	and the parameter $\mu$ is selected together with the graph by the
	constraint $\int_{\RN}|u|^2\,dx=2m^2$, which prescribes the $L^2$-size
	of the height function and not the area of the graph. The regularity of the solutions and
	the strict gradient bound recalled below justify this interpretation
	for the solutions considered here.
	
	Normalized solutions of nonlinear Schr\"odinger equations, that is, of \eqref{eq} with the Born--Infeld operator replaced by $-\Delta$, have been intensively studied since \cite{J97}. We refer to \cite{BS,BM,JL20,BMRV,BdS} for the $L^2$-supercritical case, to \cite{Soave20,Soave20crit,JJLV22,JeLe22,WW22} for combined power nonlinearities, as in Example \ref{examples of f}\ref{ex combined}, {to} \cite{JL22,CGT} for the cubic--quintic case, and to \cite{HT,CGIT,CGITGS,CGT,DoSeTi} for the Lagrangian and the action--energy points of view adopted here. For the Born--Infeld equation with prescribed $L^2$-norm, Baldelli, the first author and Pomponio \cite{BMP} obtained solutions in the $L^2$-subcritical case with a pure power nonlinearity, by truncating the operator and removing the truncation afterwards by means of uniform estimates. 
	The $L^2$-supercritical case is not covered by this argument
	and was left open in \cite{BMP}.
	In this paper we solve the $L^2$-supercritical and the $L^2$-critical problems for all masses above a threshold, which is sharp within the class of action ground states, and we establish a duality between action and energy ground states which covers the $L^2$-critical case as well. We rely on the existence and regularity theory for the prescribed frequency problem developed in \cite{BIMM,BIM}, which provides, in particular, the strict gradient bound $\lVert\n u\rVert_\infty<1$ for its critical points. For the corresponding regularity theory of the electrostatic Born--Infeld equation with a prescribed charge density we refer to \cite{BCF,Bon-Iac2,Bon-Iac3,BIMM2024}.
	
	For the nonlinearity $f$ we make the following assumptions:
	
	\begin{enumerate}
		\itemcustom{\rm (F1)}{F1} $f\in \cC^1(\R)$ and $f(0)=f'(0)=0$;
		\itemcustom{\rm (F2)}{F2} $s \mapsto \frac{f(s)}{|s|}$ is strictly increasing on $\R\setminus\{0\}$;
		\itemcustom{\rm (F3)}{F3} $\displaystyle\lim_{|s|\to\infty}\frac{f(s)}s=+\infty$;
		\itemcustom{\rm (F4)}{F4} for some $2<p<2^*=\frac{2N}{N-2}$,
		\[
		0<\liminf_{s\to0}\frac{f(s)}{|s|^{p-2}s}
		\le\limsup_{s\to0}\frac{f(s)}{|s|^{p-2}s}<\infty.
		\]
	\end{enumerate}
	
	Condition \ref{F3} is the superlinearity condition (G3) of \cite{BIM} and does not require a power-like behavior of $f$ at infinity. Condition \ref{F4} only requires two-sided power bounds near the origin: it is equivalent to the existence of $c_0,C_0,s_0>0$ such that $c_0|s|^{p-2}\le f(s)/s\le C_0|s|^{p-2}$ for $0<|s|\le s_0$, and the quotient $f(s)/(|s|^{p-2}s)$ need not have a limit as $s\to0$. The exponent $p$ decides whether the problem is $L^2$-{\em subcritical}, $L^2$-{\em critical} or $L^2$-{\em supercritical}, according as $p<2+\frac4N$, $p=2+\frac4N$ or $p>2+\frac4N$. No upper bound on the growth of $f$ at infinity
	is required: the gradient constraint and the prescribed
	mass yield a uniform $L^\infty$ bound on the admissible
	functions, see Lemma \ref{Linfty bound}. Subject to the remaining assumptions,
	nonlinearities growing faster than any power are therefore allowed, see Example \ref{examples of f}\ref{ex nonpower infinity}.
	
	\subsection{Variational framework}\label{subsec framework}
	
	We briefly recall the results of \cite{BIM} on the prescribed frequency problem:
	\begin{equation}\label{eq-mu}\tag{$BI_{\mu}$}
		\begin{cases}
			-{\rm div}\left(\dfrac{\n u}{\sqrt{1-|\n u|^2}}\right)+ \mu u = f(u), & \hbox{in }\RN,\; N\geq 3,
			\\[5mm]
			u(x)\to 0, &\hbox{as }|x|\to +\infty.
		\end{cases}
	\end{equation}
	The natural setting for the Born--Infeld operator, introduced in \cite{BIMM,BIM}, is the reflexive Banach space
	\begin{equation*}
		X := \left\{u\in H^1(\RN): \nabla u \in L^{2j_0}(\RN) \right\},\qquad j_0\in\N,\;2j_0 > N{,}
	\end{equation*}
	equipped with the norm
	\begin{equation*}
		\lVert u\rVert_X := \left(\lVert u\rVert_{H^1}^{2}+\lVert \nabla u\rVert_{2j_0}^{2}\right)^{1/2},
	\end{equation*}
	where $\lVert \cdot\rVert_{r}$ denotes the standard $L^r$-norm. Thanks to $2j_0>N$, the following embedding holds (cf. \cite[Lemma 2.15]{BIMM}):
	\begin{equation}\label{embedding}
		X\subset H^1(\R^N)\cap \cC_{0}(\RN), \qquad \cC_{0}(\RN) =\{u\in\cC(\RN):\lim\limits_{|x|\to\infty} u(x)=0\}.
	\end{equation}
	The {\em action} of \eqref{eq-mu}, denoted by $I_\mu$, splits as the sum of the convex, lower semicontinuous term
	\begin{equation*}
		\Psi(u):=
		\begin{cases}
			\displaystyle \irn \left(1-\sqrt{1-|\nabla u|^2}\right)\,dx,
			&\text{if }\lVert\nabla u\rVert_{\infty}\leq 1,\\[3mm]
			+\infty, &\text{otherwise}{,}
		\end{cases}
	\end{equation*}
	and of a $\cC^1$-term
	\begin{equation*}
		\Phi_{\mu}(u):=\frac{\mu}{2}\lVert u\rVert_2^2-\irn F(u)\,dx,
		\qquad \hbox{for every }\mu\in(0,+\infty).
	\end{equation*} 
	We denote by $\D=\{u\in X:\lVert\n u\rVert_\infty\le1\}$ the effective domain and we study the {\it energy }and the {\it action} functionals
	\begin{equation}\label{def of energy and action}
		\cE(u):=\Psi(u)-\irn F(u)\,dx,\qquad
		I_\mu(u):=\Psi(u)+\frac\mu2\lVert u\rVert_2^2-\irn F(u)\,dx,\qquad u\in\D,\ \mu>0.
	\end{equation}
	The functions of prescribed mass form the set
	\begin{equation*}
		\cS_m:=\big\{u\in\D:\lVert u\rVert_2^2=2m^2\big\}.
	\end{equation*}
	Since the energy $\cE$ of \eqref{def of energy and action} is not differentiable, we do not define the solutions of \eqref{eq} as constrained critical points of $\cE$ on $\cS_m$: a {\em solution} of \eqref{eq} is a pair $(\mu,u)\in(0,+\infty)\times\cS_m$ such that $u$ is a critical point of $I_\mu$ in the sense of Szulkin \cite{S}, see \eqref{def critical point} below, the {\em frequency} $\mu$ playing the role of a Lagrange multiplier. By the results recalled below, such a $u$ is a classical solution of the equation in \eqref{eq} with $\lVert\n u\rVert_\infty<1$.
	When the frequency is specified or clear from the context, we also call
	$u$ itself a solution.
	
	The gradient constraint and the prescribed mass bound $\|u\|_\infty$
	uniformly on $\cS_m$, and consequently $\cE$ is bounded from below on
	every mass sphere, even for a purely focusing $L^2$-supercritical
	nonlinearity. Both estimates are proved in Appendix \ref{sec estimates}.
	Unlike the focusing pure power problem for the Laplacian, the
	Born--Infeld problem therefore allows one to minimize the energy globally on $\cS_m$
	also at supercritical masses. The main results below determine when
	this infimum is attained and show that the corresponding mass threshold differs from the
	threshold for normalized action ground states.
	
	In \cite[Propositions 2.10 and 2.13]{BIMM} (see also \cite[Lemma 2.3]{BIM}), it is shown that critical points of $I_{\mu}$ in the sense of Szulkin \cite{S} are weak solutions, and in fact classical solutions with $\lVert\n u\rVert_\infty<1$, of the prescribed frequency problem \eqref{eq-mu}. Recall that $u\in\D$ is a {\em critical point} of $I_\mu$ if
	\begin{equation}\label{def critical point}
		\Psi(v)-\Psi(u)+\partial_u\Phi_{\mu}(u)[v-u]
		=\Psi(v)-\Psi(u)+\mu\irn u(v-u)\,dx-\irn f(u)(v-u)\,dx\geq0
	\end{equation}
	for every $v\in X$. 
	The Poho\v{z}aev-type (or combined {Nehari--Poho\v{z}aev}) functional related to \eqref{eq-mu} is defined for any $\mu>0$ as 
	\begin{equation*}
		M_{\mu}(u) := \Psi(u)-\irn H_{\mu}(u)\,dx,
	\end{equation*}
	where 
	\begin{equation}\label{function Hmu}
		H_{\mu}(s):= F(s)+\frac{1}{N}f(s)s-\mu\left(\frac{1}{2}+\frac{1}{N}\right) |s|^2 .
	\end{equation}
	According to \cite[Lemma 2.3]{BIM} (see also \cite[Proposition 2.6]{BIMM}), every critical point {$u$} of $I_{\mu}$ satisfies the Poho\v{z}aev identity $M_{\mu}(u)=0$. Let 
	\begin{equation*}
		\cM_{\mu} := \left\{u\in \D\setminus\{0\} : M_{\mu}(u)=0 \right\}.
	\end{equation*}
	Hence $\cM_{\mu}$ contains all nontrivial critical points of $I_{\mu}$. Similarly to \cite[Definition 1.1]{DoSeTi}, we define the {\it action ground state level} for a given $\mu$ as
	\begin{equation}\label{def of g_action}
		a(\mu):= \inf_{u\in \cM_{\mu} } I_{\mu}(u).
	\end{equation}
	Assumptions \ref{F1}--\ref{F3} coincide with the assumptions (G1)--(G3) of \cite{BIM} for the equation with unit frequency ($\mu=1$). For an arbitrary $\mu>0$, the change of variables $v(x)=\sqrt\mu\,u(x/\sqrt\mu)$ transforms \eqref{eq-mu} into the equation with unit frequency and with $f$ replaced by $g_\mu(s)=\mu^{-1/2}f(s/\sqrt\mu)$, which satisfies the same assumptions, see also {Appendix \ref{sec ground states}}.
	In particular, $a(\mu)>0$ by \cite[Lemma 5.4]{BIM}, and \cite[Theorem 1.6 and Section 5]{BIM} apply to $I_{\mu}$: $a(\mu)$ is attained by some $u_{\mu}\in\cM_{\mu}$, every minimizer is a critical point of $I_{\mu}$ in the sense of \eqref{def critical point} (see \cite[Lemma 5.5]{BIM}, or apply \cite[Lemma A.1]{MedZeng2026}), and a radial minimizer can be chosen. In Appendix \ref{sec ground states} we supplement \cite{BIM} by showing that every minimizer is radially symmetric about some point, see Theorem \ref{symmetry of action ground states}.
	
	A minimizer $u\in\cM_\mu$ satisfying $I_\mu(u)=a(\mu)$ is called
	an {\it action ground state solution} at frequency $\mu$.
	By the results just recalled, these are precisely the nontrivial critical
	points of $I_\mu$, in the sense of \eqref{def critical point}, with least action among all nontrivial critical
	points. We denote their set by
	\begin{equation*}
		\cC(\mu) := \left\{u\in \cM_{\mu} : I_{\mu}(u)=a(\mu)\right\}.
	\end{equation*}
	This set is nonempty, and every element is radial up to a translation.
	If in addition $u\in\cS_m$, we call $u$ a
	{\it normalized action ground state solution} of mass $2m^2$.
	The associated pair $(\mu,u)$ solves \eqref{eq}.
	
	\subsection{Main results}
	
	We introduce
	\begin{equation}\label{def of g_energy}
		b_m(\mu):=a(\mu)-\mu m^2.
	\end{equation}
	Following \cite{DoSeTi,CGT}, we look for solutions of \eqref{eq} among the action ground states of \eqref{eq-mu}, choosing the frequency $\mu$ in such a way that the mass is the prescribed one. Since neither the uniqueness of the action ground states nor the differentiability of $a$ is known, we consider {half of the smallest and of the largest mass} in $\cC(\mu)$,
	\begin{equation}\label{def of g_mass}
		c_-(\mu):=\inf_{u\in\cC(\mu)}
		\frac12\lVert u\rVert_2^2,
		\qquad
		c_+(\mu):=\sup_{u\in\cC(\mu)}
		\frac12\lVert u\rVert_2^2{,}
	\end{equation}
	and we prove in Theorem \ref{relation btw eng' and mass} that $a$ has one-sided derivatives, given by
	\begin{equation}\label{one-side derivatives}
		a'_+(\mu)=c_-(\mu),\qquad a'_-(\mu)=c_+(\mu).
	\end{equation}
	If $\mu_0>0$ is a local minimum point of $b_m$, these formulas give
	\[
	0\le (b_m)'_+(\mu_0)=c_-(\mu_0)-m^2
	\le c_+(\mu_0)-m^2=(b_m)'_-(\mu_0)\le0,
	\]
	hence $c_-(\mu_0)=c_+(\mu_0)=m^2$: every $u\in\cC(\mu_0)$ has mass $2m^2$, and $(\mu_0,u)$ solves \eqref{eq}. Our first main result describes the masses for which this happens.
	
	\begin{theorem}\label{main theorem A}
		Assume \ref{F1}--\ref{F4} with $2+\frac4N\le p<2^*$. {Then $c_*:=\inf_{\mu>0}c_-(\mu)>0$ and it is attained if $p>2+\frac4N$.}  For every $m^2>c_*$, and for $m^2=c_*$ whenever $c_*$ is attained, problem \eqref{eq} admits a {radial solution of constant sign with a positive Lagrange multiplier}. More precisely:
		\begin{enumerate}[label={\rm (\roman*)}, ref={\rm (\roman*)}]
			\item\label{main A i} if $m^2<c_*$, or if $m^2=c_*$ and $c_*$ is not attained, then no action ground state $u\in\cC(\mu)$, $\mu>0$, satisfies $\lVert u\rVert_2^2=2m^2$;
			\item\label{main A ii} if $m^2=c_*=c_-(\mu_0)$ for some $\mu_0>0$, {then problem \eqref{eq} has a radial solution $(\mu_{0},u)$ of constant sign with $u\in \cC(\mu_0)$. Moreover, it holds that $\cE(u)=b_m(\mu_0)\ge0$, and the inequality is strict if $p>2+\frac4N$};
			\item\label{main A iii} if $m^2>c_*$, then there is $\mu>0$ such that every $u\in\cC(\mu)$ satisfies $\frac12\lVert u\rVert_2^2=c_-(\mu)=c_+(\mu)=m^2$. In particular, \eqref{eq} has a radial solution $(\mu,u)$ of constant sign with $\cE(u)=b_m(\mu)$.
		\end{enumerate}
		Moreover, if $f$ is odd, all the solutions in \ref{main A ii} and \ref{main A iii} can be taken positive.
	\end{theorem}
	
	{Remark that whether problem \eqref{eq} has a solution in the range $m^2<c_*$ is still open.} We next present the existence result in the $L^2$-subcritical case.
	
	\begin{theorem}\label{main theorem sub}
		Assume \ref{F1}--\ref{F4} with $2<p<2+\frac4N$. Then for every $m>0$ there is $\mu>0$ such that every $u\in\cC(\mu)$ satisfies $\frac12\lVert u\rVert_2^2=c_-(\mu)=c_+(\mu)=m^2$. In particular, \eqref{eq} has a radial solution $(\mu,u)$ of constant sign, positive if $f$ is odd, with $\cE(u)=b_m(\mu)=\min_{(0,\infty)}b_m<0$.
	\end{theorem}
	
	In the $L^2$-subcritical case, Theorem \ref{main theorem sub} recovers the existence result of \cite{BMP} without truncating the operator and without any upper restriction on the growth of $f$ at infinity. In the $L^2$-supercritical case, Theorem \ref{main theorem A} provides normalized solutions for every $m^2\ge c_*$: the solutions for $m^2>c_*$ come from a local minimum point of $b_m$, while for $m^2=c_*$ a ground state of minimum mass is selected. In the $L^2$-critical case the same holds for every $m^2>c_*$, but Lemma \ref{behavior at zero} only bounds the masses of the action ground states away from zero at small frequencies, instead of forcing them to infinity. Hence the attainment of $c_*$ is not automatic, and the mass $2c_*$ is attained by an action ground state exactly when $c_*$ is attained. We stress that the nonexistence statement \ref{main A i} concerns action ground states only, and not arbitrary solutions of \eqref{eq}.
	
	The solutions of Theorems \ref{main theorem A} and \ref{main theorem sub} are action ground states of \eqref{eq-mu}, which have constant sign by Theorem \ref{symmetry of action ground states}, but for a non-odd $f$ they need not be positive: if $-f(-s)>f(s)$ for every $s>0$, then $F(-s)>F(s)$ for $s>0$, hence $\Lambda(\mu,-u)<\Lambda(\mu,u)$ for every $u\ge0$, $u\neq0$, with $\Lambda$ as in \eqref{def of max-energy map}, and all the action ground states are negative. When $f$ is odd, the solutions can be taken positive. For a general $f$, positive solutions are obtained by applying Theorems \ref{main theorem A}--\ref{main theorem energy} to the odd extension of $f|_{[0,+\infty)}$, which satisfies \ref{F1}--\ref{F4} as well. The mass thresholds in this application are those associated with the odd extension. The resulting positive solutions solve the original equation, but need not be action or energy ground states for the original non-odd nonlinearity.
	
	We next investigate the relation
	between the normalized action ground states obtained in
	Theorems \ref{main theorem A} and \ref{main theorem sub}
	and the minimization of the energy at prescribed mass.
	For $m>0$, let
	\begin{equation}\label{def of energy level}
		e_m:=\inf_{u\in\cS_m}\cE(u),
	\end{equation}
	which is finite by Lemma \ref{energy bounded below} in all three
	mass regimes, in the $L^2$-critical and in the $L^2$-supercritical case as well.
	A minimizer $u\in\cS_m$ satisfying $\cE(u)=e_m$ is called
	an {\it energy ground state} of mass $2m^2$.
	For the semilinear $L^2$-critical problem, a comparison between the minimization at fixed mass and the Lagrangian minimax levels is given in \cite[Theorem 1.6 and Remark 1.11]{CGITGS}.
	The following theorem characterizes the level $e_m$ and
	its attainment, and it identifies the normalized
	action ground states which are energy ground states.
	
	\begin{theorem}[energy ground states]\label{main theorem energy}
		Assume \ref{F1}--\ref{F4} with $2<p<2^*$ and let $m>0$.
		\begin{enumerate}[label={\rm (\roman*)}, ref={\rm (\roman*)}]
			\item\label{energy i} $e_m=\inf_{\mu>0}b_m(\mu)$, and $m^2\mapsto e_m$ is concave and nonincreasing.
			\item\label{energy ii} $e_m$ is attained if and only if $b_m$ attains its infimum on $(0,\infty)$. In this case, every global minimum point $\mu_m$ of $b_m$ satisfies $c_-(\mu_m)=c_+(\mu_m)=m^2$, and $u\in\cS_m$ is an energy ground state if and only if $u\in\cC(\mu_m)$ for some global minimum point $\mu_m$ of $b_m$. In particular, $(\mu_m,u)$ solves \eqref{eq}, and $u$ has constant sign and is radially symmetric about some point.
			\item\label{energy iii} If $2<p<2+\frac4N$, then $e_m<0$ is attained for every $m>0$. If $2+\frac4N\le p<2^*$, set
			\[
			r_e:=\inf_{\mu>0}\frac{a(\mu)}{\mu}.
			\]
			Then $c_*\le r_e<\infty$, with $c_*$ as in Theorem \ref{main theorem A}, and $c_*<r_e$ if $p>2+\frac4N$. The infimum $e_m=0$ is not attained for $m^2<r_e$, whereas $e_m<0$ is attained for $m^2>r_e$. For $m^2=r_e$, $e_m=0$ is attained if and only if $\mu\mapsto a(\mu)/\mu$ attains its infimum, which is the case if $p>2+\frac4N$.
		\end{enumerate}
	\end{theorem}
	
	Theorem \ref{main theorem energy}\ref{energy i} identifies
	$-e_m=\sup_{\mu>0}\big(\mu m^2-a(\mu)\big)$, as a function of $m^2$, with the
	{\it Legendre--Fenchel transform} of $a$. This is the Born--Infeld
	counterpart of \cite[Theorem 1.2]{DoSeTi}, while \ref{energy ii}
	corresponds to the rigidity result \cite[Theorem 1.3]{DoSeTi}.
	The rays $\{tu:t>0\}$ on which the Nehari approach of \cite{DoSeTi} is based
	leave $\D$ for large $t$, and we use the combined
	Nehari--Poho\v{z}aev constraint instead. The resulting two-threshold
	structure appears even for a single focusing power: for
	$c_*<m^2<r_e$ there are action ground states of positive
	energy but no energy ground state, whereas for $m^2>r_e$ the
	energy minimum is attained. In the $L^2$-supercritical case $c_*<r_e$ and both infima are attained. In the $L^2$-critical case our estimates decide neither whether $c_*<r_e$ nor whether the two infima are attained. A similar distinction between
	normalized solutions and energy minimizers is known for
	focusing--defocusing nonlinearities: for the Schr\"odinger equation with a focusing cubic and a defocusing quintic term on $\R^3$, the infimum of the energy on the mass sphere is attained only for masses at and above a threshold, whereas normalized solutions of positive energy, which are not energy minimizers, exist as well, see \cite{JL22,CGT}.
	In the present problem the boundedness of the energy and the
	separation of the two thresholds are caused by the gradient
	constraint alone, without any defocusing term.
	
	For the $L^2$-critical pure power, the threshold can be computed explicitly, and the nonattainment at the threshold can be decided, by means of the sharp
	Gagliardo--Nirenberg inequality of Weinstein \cite{Weinstein}.
	
	\begin{corollary}[$L^2$-critical pure power]\label{critical power corollary}
		Let $f(s)=|s|^{4/N}s$, and let $Q$ be a positive ground state of
		\[
		-\Delta Q+Q=Q^{1+4/N}\quad\text{in }\RN.
		\]
		Set $m_c:=\lVert Q\rVert_2/\sqrt2$. Then
		$r_e=m_c^2$, with $r_e$ as in Theorem
		\ref{main theorem energy}\ref{energy iii}, and
		\begin{enumerate}[label={\rm (\roman*)}, ref={\rm (\roman*)}]
			\item for $0<m\le m_c$, $e_m=0$ and the infimum is not attained;
			\item for $m>m_c$, $-\infty<e_m<0$ and the infimum is attained.
		\end{enumerate}
		In particular, $\mu\mapsto a(\mu)/\mu$ does not attain its infimum.
	\end{corollary}
	
	For comparison, if the Born--Infeld density is replaced by the Dirichlet density $|\n u|^2/2$
	and $\cS_m$ by the mass sphere in $H^1(\RN)$, without any gradient
	constraint, then the infimum of the energy is zero for $0<m\le m_c$, attained only when $m=m_c$, and
	it is $-\infty$ for $m>m_c$, see e.g. \cite{Weinstein,Caz}. Thus the Born--Infeld term changes both
	the behavior at the threshold and the existence of minimizers
	above it. Corollary \ref{critical power corollary} is proved at
	the end of Section \ref{sec structure}.

	Let us describe the strategy of the proofs. The dilations $u_t(x)=t u(x/t)$ preserve the gradient constraint and provide a projection onto the combined Nehari--Poho\v{z}aev set $\cM_\mu$. {Comparing the maxima of $I_\mu$ along these dilations at two frequencies, see Lemma \ref{comparison lemma}, we obtain
	\begin{equation*}
		\limsup_{h\to0^+}\frac{a(\mu+h)-a(\mu)}{h}\le c_-(\mu),
		\qquad
		\liminf_{h\to0^+}\frac{a(\mu)-a(\mu-h)}{h}\ge c_+(\mu)
		\qquad\text{for every }\mu>0,
	\end{equation*}
	see \eqref{upper bounds for one-sided derivatives}. These two inequalities already imply that $c_-(\mu)=c_+(\mu)=m^2$ at every local minimum point $\mu$ of $b_m$. Such a point exists in $(\mu_1,\infty)$ whenever $c_-(\mu_1)<m^2$: in this case $b_m(\mu)<b_m(\mu_1)$ for $\mu>\mu_1$ close to $\mu_1$, and $b_m(\mu)\to+\infty$ as $\mu\to\infty$, because the gradient constraint forces $a(\mu)/\mu\to\infty$, see Lemma \ref{behavior at infinity}. This gives the normalized action ground states of Theorem \ref{main theorem A}\ref{main A iii}, in the $L^2$-critical and in the $L^2$-supercritical case alike, while in the $L^2$-subcritical case $b_m$ is negative at small frequencies and its minimum is global. We also prove that the ground states of the prescribed frequency problem obtained in \cite{BIM} are compact, up to translations, when the frequency varies in a compact subset of $(0,\infty)$. This compactness gives the reverse inequalities, hence the one-sided derivative formulas \eqref{one-side derivatives}, and the lower semicontinuity of $c_-$, on which the attainment of $c_*$ for $p>2+\frac4N$ relies. The relation between the action and the energy ground state levels does not depend on these formulas. It follows from two elementary observations: if an action ground state at the frequency $\mu$ is dilated onto $\cS_m$, then its energy does not exceed $b_m(\mu)$, and conversely} every $u\in\cS_m$ with $\cE(u)\le0$ belongs to $\cM_\mu$ for a suitable frequency $\mu>0$.
	
	{The Lagrangian constructions of \cite{CGIT,CGITGS,CGT} are carried out in the radial space $H^1_r(\RN)$, whereas here $a(\mu)$ and $\cC(\mu)$ are defined in the whole space $X$, without any symmetry restriction.} In the $L^2$-supercritical construction of \cite{BMP}, the boundedness of the Palais--Smale sequences satisfying an asymptotic Poho\v{z}aev identity relies on a coercivity condition, denoted there by (b3), which the truncated Born--Infeld operators do not satisfy near the singular gradient level, see \cite[Remark 5.6]{BMP}. For this reason the supercritical results of \cite{BMP} do not pass to the Born--Infeld limit. Our approach avoids this difficulty: we work with the ground states of the untruncated prescribed frequency problem and select their mass afterwards, and no estimate at fixed mass for the truncated problems is needed.

	The following examples illustrate the scope of \ref{F1}--\ref{F4}.
	
	\begin{example}\label{examples of f}
		Unless stated otherwise, $f$ is given for $s\ge0$ and extended oddly
		to $\R$, with $f(0)=0$. In \ref{ex nonpower infinity}--\ref{ex asymmetric},
		let $2<p<2^*$.
		\begin{enumerate}[label={\rm (\roman*)}, ref={\rm (\roman*)}]
			\item\label{ex combined} {\it Pure and combined powers.}
			Let $f(s)=\sum_{i=1}^ka_is^{p_i-1}$, where $k\ge1$, $a_i>0$,
			$2<p_1<\cdots<p_k$ and $p_1<2^*$.
			Then \ref{F1}--\ref{F4} hold with $p=p_1$ and
			$f(s)/s^{p_1-1}\to a_1$ as $s\to0^+$.
			The higher powers may equal or exceed $2^*$.
			\item\label{ex nonpower infinity} {\it Non-power growth at infinity.}
			\[
			f(s)=s\log(1+s^{p-2})
			\quad\hbox{and}\quad
			f(s)=s^{p-1}e^{s^2}.
			\]
			Both satisfy \ref{F1}--\ref{F4}, with $f(s)/s^{p-1}\to1$
			as $s\to0^+$. At infinity, the first grows more slowly than
			any superlinear power, whereas the second grows faster than any power.
			\item\label{ex oscillatory zero} {\it No asymptotic coefficient at zero.}
			\[
			f(s)=s^{p-1}\big[2+\sin((p-2)\log s)\big],\qquad s>0.
			\]
			The odd extension is of class $\cC^1$ with $f'(0)=0$.
			Writing $\vartheta=(p-2)\log s$, we have
			\[
			\left(\frac{f(s)}s\right)'
			=(p-2)s^{p-3}(2+\sin\vartheta+\cos\vartheta)>0.
			\]
			Thus \ref{F1}--\ref{F4} hold, but the lower and upper limits
			in \ref{F4} are $1$ and $3$.
			\item\label{ex asymmetric} {\it Non-odd nonlinearities.}
			\[
			f(s)=
			\begin{cases}
				a_+s^{p-1},&s\ge0,\\
				-a_-|s|^{p-1},&s<0,
			\end{cases}
			\qquad a_+,a_->0,\quad a_+\ne a_-.
			\]
			Here \ref{F1}--\ref{F4} hold with lower and upper limits
			$\min\{a_+,a_-\}$ and $\max\{a_+,a_-\}$ in \ref{F4}.
		\end{enumerate}
	\end{example}

	The symmetry of every action ground state, rather than merely the existence of a radial one, is established in Appendix \ref{sec ground states}. The proof runs as follows. The minimality implies that {for an action ground state $u$ of constant sign $\varepsilon\in\{-1,1\}$, the rearrangement $\varepsilon|u|^*$ is again an action ground state with the same kinetic energy}, and since the Born--Infeld density is a series of even powers of $|\nabla u|$ with positive coefficients, this reduces to the equality case of the P\'olya--Szeg\H{o} inequality for the Dirichlet integral, settled by Brothers and Ziemer \cite{BZ}, which we use in the quantitative form of Burchard and Ferone \cite{BF}. A uniqueness argument for the radial ordinary differential equation then gives the strict monotonicity of the radial profile, see Theorem \ref{symmetry of action ground states}. 
	
	The paper is organized as follows. Section \ref{sec action} proves the one-sided differentiability of $a(\mu)$. Section \ref{sec structure} derives from the gradient constraint the asymptotics of $a(\mu)$ and $c_\pm(\mu)$ and contains the proofs of Theorems \ref{main theorem A}--\ref{main theorem energy}. Appendix \ref{sec estimates} contains the two elementary estimates \eqref{Linfty by L2} and \eqref{lower bound of energy} on which Sections \ref{sec action}--\ref{sec structure} rely. Appendix \ref{sec ground states} proves the sign, the symmetry and the monotonicity of every Born--Infeld action ground state. Since this result applies to every ground state, we can center by a translation the ground states selected by their mass, and we obtain the symmetry statements for all the energy ground states.
	
	\smallskip
	\noindent{\bf Notation.} $\lVert\cdot\rVert_r$ is the norm of $L^r(\RN)$, $B_R(x)$ is the open ball of radius $R$ centered at $x$, $B_R=B_R(0)$, $|E|$ is the Lebesgue measure of $E$, {$(v)_+=\max\{v,0\}$ is the positive part of $v$}, $X_r$ is the subspace of radial functions in $X$, and $C$ denotes positive constants which may change from line to line.
	\section{Action ground state level}\label{sec action}
	
	In this section, we study the action ground state level as a function of the frequency.
	The following theorem is the Born--Infeld counterpart of \cite[Theorem 1.3]{CGT} and of \cite[Theorem 1.5]{DoSeTi}. It is the first main ingredient of our approach.
	Throughout this section we assume \ref{F1}--\ref{F3}.
	\begin{theorem}\label{relation btw eng' and mass}
		Let $a$, $c_-$ and $c_+$ be the functions defined in \eqref{def of g_action} and \eqref{def of g_mass}. Then
		\begin{equation*}
			a_{+}'(\mu) = c_{-}(\mu), \qquad  a_{-}'(\mu) = c_{+}(\mu)
		\end{equation*}
		for $\mu \in (0, \infty)$, where $a_{+}'(\mu)$ and $a_{-}'(\mu)$ are the right and left derivatives, respectively.
	\end{theorem}
	
	To prove Theorem \ref{relation btw eng' and mass}, we need some preliminaries, most of which are taken from \cite{BIM}.
	
	Let us point out that no upper bound on the growth of $f$ at infinity is needed in what follows. Indeed, by \ref{F1}, for every $R>0$ and every compact interval $J\subset(0,\infty)$ there is a constant $C_{R,J}$ such that
	\[
	|F(s)|+|f(s)s|+|H_\mu(s)|\le C_{R,J}s^2,
	\qquad |H_\mu'(s)|\le C_{R,J}|s|
	\]
	for $|s|\le R$ and $\mu\in J$. Since $X\hookrightarrow L^\infty(\RN)$, these local estimates justify all the nonlinear integrals below and their differentiation on bounded subsets of $X$.
	
	For $u\in\D$ and $\t>0$, set $u_{\t}(x):=\t u(x/\t)$. This scaling preserves the gradient constraint and satisfies
	\begin{equation*}
		\Psi(u_{\t})=\t^N\Psi(u),
		\qquad \lVert u_{\t}\rVert_2^2=\t^{N+2}\lVert u\rVert_2^2.
	\end{equation*}
	{Hence, for every $\mu>0$, $u\in\D\setminus\{0\}$ and $\t>0$,
	\begin{equation}\label{fiber map of action}
		I_{\mu}(u_{\t})=\t^N\left(\Psi(u)+\frac{\mu\t^2}{2}\lVert u\rVert_2^2-\irn F(\t u)\,dx\right),
		\qquad
		\frac{d}{d\t}I_{\mu}(u_{\t})=N\t^{N-1}\vp(\t,\mu,u),
	\end{equation}
	where
	\begin{equation*}
		\vp(\t,\mu,u):=\t^{-N}M_{\mu}(u_{\t})
		=\Psi(u)-\irn H_{\mu}(\t u)\,dx.
	\end{equation*}}%
	In view of \ref{F1}--\ref{F3} and \cite[Lemma 5.2]{BIM}, for every fixed $\mu>0$ and $u\in\D\setminus\{0\}$ there exists a unique $\t_{\mu}(u)>0$ such that $u_{\t_{\mu}(u)}\in\cM_{\mu}$. {Moreover,
	\begin{equation}\label{sign of fiber derivative}
		\vp(\t,\mu,u)>0\quad\text{for }0<\t<\t_{\mu}(u),
		\qquad
		\vp(\t,\mu,u)<0\quad\text{for }\t>\t_{\mu}(u).
	\end{equation}
	Indeed, $\vp(\t_{\mu}(u),\mu,u)=0$ and, in view of \eqref{function Hmu}, the map $\t\mapsto\t^{-2}\vp(\t,\mu,u)$ is strictly decreasing, because $\Psi(u)>0$ and, by \ref{F2}, the functions $F(s)/s^2$ and $f(s)/s$ are nondecreasing in $|s|$ on each of the half-lines $(-\infty,0)$ and $(0,\infty)$. By \eqref{fiber map of action} and \eqref{sign of fiber derivative}, the map $\t\mapsto I_{\mu}(u_{\t})$ is strictly increasing on $(0,\t_{\mu}(u)]$ and strictly decreasing on $[\t_{\mu}(u),\infty)$, and it is positive on $(0,\t_{\mu}(u)]$, because $I_{\mu}(u_{\t})\to0$ as $\t\to0^+$. In particular, $\t_{\mu}(u)=1$ for $u\in\cM_{\mu}$ and
	\begin{equation}\label{strict maximum on fibers}
		I_{\mu}(u_{\t})<I_{\mu}(u)
		\qquad\text{for every }u\in\cM_{\mu}\text{ and }\t\in(0,\infty)\setminus\{1\}.
	\end{equation}}%
	We define the projection
	\begin{equation*}
		m_{\mu}(u) := u_{\t_{\mu}(u)},
		\qquad u\in\D\setminus\{0\},
	\end{equation*}
	and the maximal-energy map
	\begin{equation}\label{def of max-energy map}
		\Lambda(\mu,u):= I_{\mu}(m_{\mu}(u))= I_{\mu}(u_{\t_{\mu}(u)}){{}=\max_{\t>0}I_{\mu}(u_{\t})}.
	\end{equation}
	{Thus $a(\mu)=\inf_{u\in\D\setminus\{0\}}\Lambda(\mu,u)$, and $\Lambda(\mu,u)= a(\mu)$ if $u\in \cC(\mu)$. The next lemma shows that $\t_{\mu}(u)$ depends continuously on $\mu$ and $u$, and it compares the action ground state levels at two frequencies.}
	
	\begin{lemma}\label{comparison lemma}
		{Let $\mu_n\to\mu_0\in(0,\infty)$ and let $w_n,w_0\in\D\setminus\{0\}$ be such that $w_n\to w_0$ in $X$. Then $\t_{\mu_n}(w_n)\to\t_{\mu_0}(w_0)$. Moreover,
		\begin{equation}\label{comparison of levels}
			a(\mu)\le a(\nu)+\frac{\mu-\nu}{2}\,\t_{\mu}(u)^{N+2}\lVert u\rVert_2^2
		\end{equation}
		for every $\mu,\nu>0$ and $u\in\cC(\nu)$.}
	\end{lemma}
	\begin{proof}
		{Fix $\t>0$. Since $\Psi(w_n)\to\Psi(w_0)$ by \cite[Lemma 2.1]{BIM}, and since $w_n\to w_0$ in $L^2(\RN)$ and $(w_n)$ is bounded in $L^\infty(\RN)$, the local estimates above yield
		\begin{equation*}
			\vp(\t,\mu_n,w_n)\to\vp(\t,\mu_0,w_0)
			\qquad\text{as }n\to\infty.
		\end{equation*}
		Let $\t_0:=\t_{\mu_0}(w_0)$ and $0<\eps<\t_0$. By \eqref{sign of fiber derivative}, we have $\vp(\t_0-\eps,\mu_0,w_0)>0>\vp(\t_0+\eps,\mu_0,w_0)$, and hence $\vp(\t_0-\eps,\mu_n,w_n)>0>\vp(\t_0+\eps,\mu_n,w_n)$ for $n$ large. Using \eqref{sign of fiber derivative} again, we get $\t_0-\eps<\t_{\mu_n}(w_n)<\t_0+\eps$ for $n$ large, which proves the first assertion.

		Let now $\mu,\nu>0$, $u\in\cC(\nu)$ and $\t:=\t_{\mu}(u)$. Since $u_{\t}\in\cM_{\mu}$ and $u\in\cM_{\nu}$, \eqref{strict maximum on fibers} gives
		\begin{equation*}
			a(\mu)\le I_{\mu}(u_{\t})
			=I_{\nu}(u_{\t})+\frac{\mu-\nu}{2}\lVert u_{\t}\rVert_2^2
			\le I_{\nu}(u)+\frac{\mu-\nu}{2}\,\t^{N+2}\lVert u\rVert_2^2,
		\end{equation*}
		which is \eqref{comparison of levels}, because $I_{\nu}(u)=a(\nu)$.}
	\end{proof}
	
	\begin{lemma}\label{Uniform local bounds for ground-state levels}
		Let $\mu_n\to\mu_0>0$. Then {$\limsup_{n\to\infty}a(\mu_n)\le a(\mu_0)$, and $(a(\mu_n))$ is bounded away from zero}.
	\end{lemma}

	\begin{proof}
		{Fix $u_0\in\cC(\mu_0)$. By Lemma \ref{comparison lemma}, we have $\t_{\mu_n}(u_0)\to\t_{\mu_0}(u_0)=1$, and \eqref{comparison of levels} with $\nu=\mu_0$ gives
		\begin{equation*}
			a(\mu_n)\le a(\mu_0)+\frac{\mu_n-\mu_0}{2}\,\t_{\mu_n}(u_0)^{N+2}\lVert u_0\rVert_2^2.
		\end{equation*}
		Letting $n\to\infty$, we get the first assertion.}

		We now prove the lower bound. {Choose $\mu_-,\mu_+>0$ such that $\mu_n\in[\mu_-,\mu_+]$ for all sufficiently large $n$, and set $G(s):=F(s)+\frac1Nf(s)s$.} In view of \ref{F1}, there exists $\delta>0$ such that
		\begin{equation}\label{bound G}
			G(s)\le \frac{(N+2)\mu_-}{4N}s^2
			\qquad\text{for }|s|\le\delta.
		\end{equation}
		We claim that every $u\in\cM_\mu$, with $\mu\in[\mu_-,\mu_+]$, satisfies
		\begin{equation}\label{lower uniform boundedness for u}
			\lVert u\rVert_\infty>\delta.
		\end{equation}
		Otherwise, by \eqref{bound G}, one has
		\begin{equation*}
			\begin{aligned}
				M_\mu(u)
				&=\Psi(u)+\frac{N+2}{2N}\mu\lVert u\rVert_2^2-\irn G(u)\,dx\\
				&\ge\Psi(u)+\frac{(N+2)\mu_-}{4N}\lVert u\rVert_2^2>0,
			\end{aligned}
		\end{equation*}
		which contradicts $u\in\cM_\mu$.
		
		For $u\in\D$, we have $\lVert\nabla u\rVert_\infty\le1$ and
		\begin{equation}\label{lowe bound of Psi}
			\Psi(u)\ge\frac12\lVert\nabla u\rVert_2^2,
		\end{equation}
		since $1-\sqrt{1-t^2}\ge\frac12t^2$ for $0\le t\le1$. Combining \eqref{lower uniform boundedness for u}, \eqref{lowe bound of Psi} and the second inequality in \eqref{Linfty by L2}, we get
		\begin{equation*}
			\Psi(u)\ge\eta>0,
		\end{equation*}
		where $\eta$ is independent of $\mu\in[\mu_-,\mu_+]$.
		
		Moreover, \ref{F2} implies
		\begin{equation}\label{corrollary of F2}
			\frac12f(s)s-F(s)\ge0
			\qquad\text{for every }s\in\R.
		\end{equation}
		For every $u\in\cM_\mu$, the identity $M_\mu(u)=0$ yields
		\begin{equation*}
			I_\mu(u)=\frac{2}{N+2}\left[\Psi(u)+\irn\left(\frac12f(u)u-F(u)\right)\,dx\right]
			\ge\frac{2\eta}{N+2}.
		\end{equation*}
		Taking the infimum over $u\in\cM_\mu$, we obtain the lower bound.
	\end{proof}

	\begin{lemma} \label{compactness lemma}
		If $\mu_{n}\to \mu_{0}\in (0,\infty)$ and $u_n \in \cC(\mu_{n})$ for every $n\in \N$, then, after passing to a subsequence, there exist $\{y_n\}_{n\in \N}\subset \RN$ and $u_0\in \cC(\mu_0)$ such that 
		\begin{equation*}
			u_{n}(\cdot + y_n) \to u_0 \qquad 
			\text{strongly in } X.
		\end{equation*}
	\end{lemma}
	\begin{proof}
		Set
		\begin{equation*}
			K(s):=\frac12f(s)s-F(s),
			\qquad G(s):=F(s)+\frac1Nf(s)s.
		\end{equation*}
		Since $u_n\in\cC(\mu_n)$, the identity $M_{\mu_n}(u_n)=0$ and \eqref{corrollary of F2} give
		\begin{equation}\label{energy identity on ground states}
			\begin{aligned}
				a(\mu_n)=I_{\mu_n}(u_n)
				&=\frac{2}{N+2}\left[\Psi(u_n)+\irn K(u_n)\,dx\right]\\
				&\ge\frac{2}{N+2}\Psi(u_n).
			\end{aligned}
		\end{equation}
		Lemma \ref{Uniform local bounds for ground-state levels} therefore yields $\Psi(u_n)\le C$. Since
		\begin{equation*}
			1-\sqrt{1-t^2}=\sum_{j=1}^{\infty}a_jt^{2j},
			\qquad a_j>0,
		\end{equation*}
		the sequences $(\lVert\nabla u_n\rVert_2)$ and $(\lVert\nabla u_n\rVert_{2j_0})$ are bounded. The Sobolev and Morrey embeddings then imply
		\begin{equation*}
			\lVert u_n\rVert_{2^*}+\lVert u_n\rVert_\infty\le C.
		\end{equation*}
		Let $R:=\sup_n\lVert u_n\rVert_\infty$. By \ref{F1}, for every $\eps>0$ there exists $C_\eps>0$ such that
		\begin{equation*}
			F(s)\le\eps s^2+C_\eps|s|^{2^*}
			\qquad\text{for }|s|\le R.
		\end{equation*}
		Choose $\mu_->0$ such that $\mu_n\ge\mu_-$ for all sufficiently large $n$, and take $\eps=\mu_-/4$. Since $I_{\mu_n}(u_n)=a(\mu_n)$, we obtain
		\begin{equation*}
			\frac{\mu_-}{4}\lVert u_n\rVert_2^2
			\le a(\mu_n)+C\lVert u_n\rVert_{2^*}^{2^*}\le C.
		\end{equation*}
		Thus, $(u_n)$ is bounded in $X$.
		
		Set
		\begin{equation*}
			\alpha_0:=\frac{N+2}{2N}\mu_0.
		\end{equation*}
		Applying the profile decomposition of \cite[Proposition 3.3]{BIM} to a subsequence, we find $(v^j)_{j\ge0}\subset X$ and translations $(y_n^j)_{n\ge1}\subset\RN$ such that
		\begin{equation}\label{profile decomposition}
			\begin{aligned}
				u_n(\cdot+y_n^j)&\weakto v^j\quad\text{in }X,
				&u_n(x+y_n^j)&\to v^j(x)\quad\text{a.e. in }\RN,\\
				\liminf_{n\to\infty}\left[\Psi(u_n)+\alpha_0\lVert u_n\rVert_2^2\right]
				&\ge\sum_{j\ge0}\left[\Psi(v^j)+\alpha_0\lVert v^j\rVert_2^2\right],\\
				\lim_{n\to\infty}\irn F(u_n)\,dx
				&=\sum_{j\ge0}\irn F(v^j)\,dx,\\
				\lim_{n\to\infty}\irn f(u_n)u_n\,dx
				&=\sum_{j\ge0}\irn f(v^j)v^j\,dx.
			\end{aligned}
		\end{equation}
		Every profile $v^j$ belongs to $\D$, since $\D$ is convex and closed in $X$, hence weakly closed.
		The same profile decomposition also gives
		\[
		\liminf_{n\to\infty}\Psi(u_n)\ge\sum_{j\ge0}\Psi(v^j).
		\]
		Indeed, the inequality in \eqref{profile decomposition} holds, with the same profiles, for every coefficient $\eps>0$ in place of $\alpha_0$, because the profiles are determined by the sequence $(u_n)$ and the proof of \cite[Proposition 3.3]{BIM} does not use the value of the coefficient. Since $\lVert u_n\rVert_2^2\le C$,
		\[
		\liminf_{n\to\infty}\Psi(u_n)+\eps C
		\ge\liminf_{n\to\infty}\big[\Psi(u_n)+\eps\lVert u_n\rVert_2^2\big]
		\ge\sum_{j\ge0}\big[\Psi(v^j)+\eps\lVert v^j\rVert_2^2\big]
		\ge\sum_{j\ge0}\Psi(v^j),
		\]
		and letting $\eps\to0^+$ proves the assertion.
		Since $(u_n)$ is bounded in $L^2$ and $\mu_n\to\mu_0$, we have
		\begin{equation*}
			M_{\mu_0}(u_n)=M_{\mu_n}(u_n)
			+\frac{N+2}{2N}(\mu_0-\mu_n)\lVert u_n\rVert_2^2\to0.
		\end{equation*}
		Combining this fact with \eqref{profile decomposition}, we obtain
		\begin{equation}\label{negative of Nehari Pohozaev}
			0\ge\sum_{j\ge0}M_{\mu_0}(v^j).
		\end{equation}
		
		Let $\mathcal J:=\{j:v^j\neq0\}$. We claim that $\mathcal J\neq\emptyset$. Indeed, if every profile vanished, then the last two identities in \eqref{profile decomposition} would imply
		\begin{equation*}
			\irn F(u_n)\,dx\to0,
			\qquad \irn f(u_n)u_n\,dx\to0.
		\end{equation*}
		Since $M_{\mu_0}(u_n)\to0$, it would follow that
		\begin{equation*}
			\Psi(u_n)+\alpha_0\lVert u_n\rVert_2^2\to0,
		\end{equation*}
		and then \eqref{energy identity on ground states} would give $a(\mu_n)\to0$. This contradicts Lemma \ref{Uniform local bounds for ground-state levels}.
		
		By \eqref{negative of Nehari Pohozaev}, there exists $j_1\in\mathcal J$ such that $M_{\mu_0}(v^{j_1})\le0$. Put
		\begin{equation*}
			\t_0:=\t_{\mu_0}(v^{j_1}).
		\end{equation*}
		{Since $\vp(1,\mu_0,v^{j_1})=M_{\mu_0}(v^{j_1})\le0$, \eqref{sign of fiber derivative} shows that $0<\t_0\le1$.} Moreover, \ref{F2} implies
		\begin{equation*}
			K(s)\ge0,
			\qquad K(\t s)\le K(s)
			\quad\text{for }s\in\R\text{ and }0<\t\le1.
		\end{equation*}
		Define
		\begin{equation*}
			\cQ(u):=\frac{2}{N+2}\left[\Psi(u)+\irn K(u)\,dx\right].
		\end{equation*}
		For $u\in\cM_\mu$, one has $I_\mu(u)=\cQ(u)$. Consequently,
		\begin{equation}\label{upper estimate for energy}
			\begin{aligned}
				a(\mu_0)
				&\le I_{\mu_0}((v^{j_1})_{\t_0})\\
				&=\frac{2\t_0^N}{N+2}
				\left[\Psi(v^{j_1})+\irn K(\t_0v^{j_1})\,dx\right]
				\le\cQ(v^{j_1}).
			\end{aligned}
		\end{equation}
		
		After passing to a further subsequence, Lemma \ref{Uniform local bounds for ground-state levels} allows us to assume that $a(\mu_n)\to L$. From \eqref{profile decomposition} and \eqref{energy identity on ground states}, we get
		\begin{equation}\label{lower bound for limit}
			L\ge\sum_{j\in\mathcal J}\cQ(v^j)
			\ge\cQ(v^{j_1})\ge a(\mu_0).
		\end{equation}
		{On the other hand, $L\le a(\mu_0)$ by Lemma \ref{Uniform local bounds for ground-state levels}.}
		Therefore, equality holds throughout \eqref{lower bound for limit}. Equality in \eqref{upper estimate for energy} forces $\t_0=1$, because $\cQ(v^{j_1})>0$ and the last inequality in \eqref{upper estimate for energy} is strict when $\t_0<1$. It follows that
		\begin{equation*}
			v^{j_1}\in\cM_{\mu_0},
			\qquad I_{\mu_0}(v^{j_1})=a(\mu_0),
		\end{equation*}
		and hence $v^{j_1}\in\cC(\mu_0)$. Furthermore, $\cQ(v)>0$ for every nonzero $v\in X$. Thus equality in \eqref{lower bound for limit} also shows that
		\begin{equation*}
			\#\mathcal J=1.
		\end{equation*}
		
		Let $u_0:=v^{j_1}$, $y_n:=y_n^{j_1}$ and $w_n:=u_n(\cdot+y_n)$. Since there is exactly one nonzero profile, \eqref{profile decomposition} gives
		\begin{equation*}
			\irn G(w_n)\,dx\to\irn G(u_0)\,dx.
		\end{equation*}
		Using $M_{\mu_n}(w_n)=0$, $\mu_n\to\mu_0$ and $M_{\mu_0}(u_0)=0$, we conclude that
		\begin{equation*}
			\Psi(w_n)+\alpha_0\lVert w_n\rVert_2^2
			\to\Psi(u_0)+\alpha_0\lVert u_0\rVert_2^2.
		\end{equation*}
		Together with $w_n\weakto u_0$ in $X$, \cite[Lemma 2.2]{BIM} now yields $w_n\to u_0$ strongly in $X$.
	\end{proof}
	
	\begin{corollary}\label{attainment of c pm}
		For every $\mu>0$ the infimum and the supremum in \eqref{def of g_mass} are attained: there exist $u_+,u_-\in\cC(\mu)$ with $\frac12\lVert u_+\rVert_2^2=c_+(\mu)$ and $\frac12\lVert u_-\rVert_2^2=c_-(\mu)$. Moreover, $\mu\mapsto c_-(\mu)$ is lower semicontinuous and $\mu\mapsto c_+(\mu)$ is upper semicontinuous on $(0,\infty)$.
	\end{corollary}
	\begin{proof}
		Let $(u_n)\subset\cC(\mu)$ be such that $\frac12\lVert u_n\rVert_2^2\to c_-(\mu)$. By Lemma \ref{compactness lemma} with $\mu_n\equiv\mu$, a subsequence of $u_n(\cdot+y_n)$ converges in $X$ to some $u_-\in\cC(\mu)$, and the translation invariance of the $L^2$-norm gives $\frac12\lVert u_-\rVert_2^2=c_-(\mu)$. In particular, $c_-(\mu)>0$, since $0\notin\cC(\mu)$. The same argument, applied to a sequence $(u_n)\subset\cC(\mu)$ with $\frac12\lVert u_n\rVert_2^2\to c_+(\mu)$, shows that $c_+(\mu)<\infty$ and that it is attained. Now let $\mu_n\to\mu_0\in(0,\infty)$ be such that $c_-(\mu_n)\to\liminf_{\mu\to\mu_0}c_-(\mu)$, and choose $u_n\in\cC(\mu_n)$ with $\frac12\lVert u_n\rVert_2^2=c_-(\mu_n)$. By Lemma \ref{compactness lemma}, a subsequence of $u_n(\cdot+y_n)$ converges in $X$ to some $u_0\in\cC(\mu_0)$, hence $\liminf_{\mu\to\mu_0}c_-(\mu)=\lim_n\frac12\lVert u_n\rVert_2^2=\frac12\lVert u_0\rVert_2^2\ge c_-(\mu_0)$. The upper semicontinuity of $c_+$ is proved in the same way.
	\end{proof}
	
	\begin{proof}[Proof of Theorem \ref{relation btw eng' and mass}]
		{Fix $\mu_0>0$. By Corollary \ref{attainment of c pm} we may choose $u_-,u_+\in\cC(\mu_0)$ such that
		\begin{equation*}
			\frac12\lVert u_-\rVert_2^2=c_-(\mu_0),
			\qquad
			\frac12\lVert u_+\rVert_2^2=c_+(\mu_0).
		\end{equation*}
		For $0<h<\mu_0$, inequality \eqref{comparison of levels} with $\nu=\mu_0$, applied to $(\mu,u)=(\mu_0+h,u_-)$ and to $(\mu,u)=(\mu_0-h,u_+)$, gives
		\begin{equation*}
			\frac{a(\mu_0+h)-a(\mu_0)}{h}\le\frac12\,\t_{\mu_0+h}(u_-)^{N+2}\lVert u_-\rVert_2^2,
			\qquad
			\frac{a(\mu_0)-a(\mu_0-h)}{h}\ge\frac12\,\t_{\mu_0-h}(u_+)^{N+2}\lVert u_+\rVert_2^2.
		\end{equation*}
		Since $\t_{\mu_0+h}(u_-)\to\t_{\mu_0}(u_-)=1$ and $\t_{\mu_0-h}(u_+)\to\t_{\mu_0}(u_+)=1$ as $h\to0^+$ by Lemma \ref{comparison lemma}, we obtain
		\begin{equation}\label{upper bounds for one-sided derivatives}
			\limsup_{h\to 0^+}\frac{a(\mu_0+h)-a(\mu_0)}{h}\le c_-(\mu_0),
			\qquad
			\liminf_{h\to0^+}\frac{a(\mu_0)-a(\mu_0-h)}{h}\ge c_+(\mu_0).
		\end{equation}

		For the reverse inequalities, let $(h_n)\subset(0,\mu_0)$ be an arbitrary sequence such that $h_n\to0$, and choose $u_n\in\cC(\mu_0+h_n)$ and $w_n\in\cC(\mu_0-h_n)$. Inequality \eqref{comparison of levels} with $\mu=\mu_0$, applied to $(\nu,u)=(\mu_0+h_n,u_n)$ and to $(\nu,u)=(\mu_0-h_n,w_n)$, gives
		\begin{equation}\label{lower bounds for one-sided derivatives}
			\frac{a(\mu_0+h_n)-a(\mu_0)}{h_n}\ge\frac12\,\t_{\mu_0}(u_n)^{N+2}\lVert u_n\rVert_2^2,
			\qquad
			\frac{a(\mu_0)-a(\mu_0-h_n)}{h_n}\le\frac12\,\t_{\mu_0}(w_n)^{N+2}\lVert w_n\rVert_2^2.
		\end{equation}
		By Lemma \ref{compactness lemma}, after passing to a subsequence, there exist $(y_n),(z_n)\subset\RN$ and $u_*,w_*\in\cC(\mu_0)$ such that
		\begin{equation*}
			u_n(\cdot+y_n)\to u_*,
			\qquad
			w_n(\cdot+z_n)\to w_*
			\qquad\text{strongly in }X.
		\end{equation*}
		Since $\t_{\mu_0}$ and the $L^2$-norm are invariant under translations, Lemma \ref{comparison lemma} yields $\t_{\mu_0}(u_n)\to\t_{\mu_0}(u_*)=1$ and $\t_{\mu_0}(w_n)\to\t_{\mu_0}(w_*)=1$. Hence, along this subsequence, the right-hand sides of \eqref{lower bounds for one-sided derivatives} converge to $\frac12\lVert u_*\rVert_2^2\ge c_-(\mu_0)$ and to $\frac12\lVert w_*\rVert_2^2\le c_+(\mu_0)$, respectively. Since every such sequence $(h_n)$ has a subsequence with these properties, we conclude that
		\begin{equation*}
			\liminf_{h\to 0^+}\frac{a(\mu_0+h)-a(\mu_0)}{h}\ge c_-(\mu_0),
			\qquad
			\limsup_{h\to0^+}\frac{a(\mu_0)-a(\mu_0-h)}{h}\le c_+(\mu_0).
		\end{equation*}
		Together with \eqref{upper bounds for one-sided derivatives}, this proves $a_+'(\mu_0)=c_-(\mu_0)$ and $a_-'(\mu_0)=c_+(\mu_0)$.}
	\end{proof}

	\section{\texorpdfstring{{Proofs of the main results}}{Proofs of the main results}}\label{sec structure}
	
	In this section we first derive the behavior of $a(\mu)$ and $c_\pm(\mu)$ as $\mu\to\infty$ and as $\mu\to0^+$, and then we prove the main results.
	
	\subsection{Behavior of $a(\mu)$, $c_\pm(\mu)$ and $b_m(\mu)$}
	
	\begin{lemma}\label{behavior at infinity}
		Assume \ref{F1}--\ref{F3} and set
		\[
		\kappa(0):=0,\qquad
		\kappa(S):=\max\left\{\frac{f(S)}S,\frac{f(-S)}{-S}\right\}
		\quad(S>0),\qquad \sigma:=\kappa^{-1}.
		\]
		Then $\kappa$ is a continuous, strictly increasing bijection of $[0,\infty)$ onto itself and, with $d_N:=|B_1|/2^{N+3}$, every $\mu>0$ and {$u\in\cM_\mu$} satisfy
		\begin{equation}\label{lower bounds at infinity}
			\lVert u\rVert_\infty\ge\sigma(\mu),\qquad
			\frac12\lVert u\rVert_2^2\ge d_N\sigma(\mu)^{N+2},\qquad
			\frac{a(\mu)}\mu\ge\frac{d_N}{2}\sigma(\mu/2)^{N+2}.
		\end{equation}
		In particular, $c_-(\mu)\to\infty$ and $a(\mu)/\mu\to\infty$ as $\mu\to\infty$.
	\end{lemma}
	\begin{proof}
		The stated properties of $\kappa$ follow from \ref{F1}--\ref{F3}: $\kappa$ is continuous and strictly increasing by \ref{F1}--\ref{F2}, $\kappa(S)\to0$ as $S\to0^+$ by \ref{F1}, and $\kappa(S)\to\infty$ as $S\to\infty$ by \ref{F3}. Let {$u\in\cM_\mu$} and set $A:=\lVert u\rVert_\infty$. Recall that $G(s)=F(s)+\frac1Nf(s)s$ and that $F(s)\le\frac12f(s)s$ by \ref{F2}, hence $G(s)\le\frac{N+2}{2N}\kappa(A)s^2$ for $|s|\le A$. Since $M_\mu(u)=0$ and $\Psi(u)\ge0$, we obtain
		\[
		\frac{N+2}{2N}\mu\lVert u\rVert_2^2
		\le\Psi(u)+\frac{N+2}{2N}\mu\lVert u\rVert_2^2
		=\irn G(u)\,dx
		\le\frac{N+2}{2N}\kappa(A)\lVert u\rVert_2^2.
		\]
		Since $u\ne0$, this gives $\kappa(A)\ge\mu$, hence $A\ge\sigma(\mu)$. Estimate \eqref{Linfty by L2} now yields
		\[
		\frac12\lVert u\rVert_2^2
		\ge\frac{|B_1|}{2^{N+3}}A^{N+2}
		\ge d_N\sigma(\mu)^{N+2}.
		\]
		{Since $\cC(\mu)\subset\cM_\mu$,} we obtain the same lower bound for $c_-(\mu)$.
		
		{Let now $w\in\cC(\mu)$ and $v:=m_{\mu/2}(w)\in\cM_{\mu/2}$, so that $\lVert v\rVert_2^2=\t_{\mu/2}(w)^{N+2}\lVert w\rVert_2^2$. Inequality \eqref{comparison of levels}, with $\mu/2$ and $\mu$ in place of $\mu$ and $\nu$, and the second inequality in \eqref{lower bounds at infinity}, applied to $v$, give
		\begin{equation*}
			a(\mu)\ge a(\mu/2)+\frac{\mu}{4}\lVert v\rVert_2^2
			\ge a(\mu/2)+\frac{d_N\mu}{2}\sigma(\mu/2)^{N+2}.
		\end{equation*}}%
		Since $a(\mu/2)>0$, the last inequality in \eqref{lower bounds at infinity} follows. Finally, $\sigma(\mu)\to\infty$ as $\mu\to\infty$, which gives both limits.
	\end{proof}
	
	We set
	\begin{equation*}
		\b:=\frac{2}{p-2}-\frac N2,
	\end{equation*}
	and we note that $\b+1>0$ if and only if $p<2^*$, and that $\b<0$ if and only if $p>2+\frac4N$.
	
	\begin{lemma}\label{behavior at zero}
		Assume \ref{F1}--\ref{F4}. Then there exist $C,c>0$ such that for sufficiently small $\mu>0$ and every $u\in\cC(\mu)$
		\begin{equation}\label{asymptotics at zero}
			a(\mu)\le C\mu^{\b+1},
		\end{equation}
		and, if $p\ge2+\frac4N$,
		\begin{equation}\label{mass at zero}
			\frac12\lVert u\rVert_2^2\ge c\mu^{\b},\qquad a(\mu)\ge c\mu^{\b+1}.
		\end{equation}
		Consequently, 
		\begin{enumerate}[label={\rm (\roman*)}, ref={\rm (\roman*)}]
			\item $a(\mu)\to0$ as $\mu\to0^+$;
			\item  if $2+\frac4N<p<2^*$, then $a(\mu)/\mu\to\infty$ and $c_-(\mu)\to\infty$ as $\mu\to0^+$;
			\item if $2<p<2+\frac4N$, then $a(\mu)/\mu\to0$ as $\mu\to0^+$.
		\end{enumerate}
		
	\end{lemma}
	\begin{proof}
		{\it Upper bound for $a(\mu)$.} Fix a nonzero function $\omega\in\cC_c^\infty(\RN)$ with $\omega\ge0$. For $\mu>0$ put
		\begin{equation*}
			w_\mu(x):=\mu^{\frac1{p-2}}\omega(\mu^{1/2}x).
		\end{equation*}
		Then $\lVert\n w_\mu\rVert_\infty=\mu^{\frac1{p-2}+\frac12}\lVert\n\omega\rVert_\infty\le1$ for small $\mu$, so $w_\mu\in\D$, and a direct computation gives
		\begin{equation}\label{scaling of trial function}
			\lVert\n w_\mu\rVert_2^2=\mu^{\b+1}\lVert\n\omega\rVert_2^2,\quad
			\lVert w_\mu\rVert_2^2=\mu^{\b}\lVert\omega\rVert_2^2,\quad
			\irn w_\mu^p\,dx=\mu^{\b+1}\irn\omega^p\,dx,\quad
			\lVert w_\mu\rVert_\infty=\mu^{\frac1{p-2}}\lVert\omega\rVert_\infty .
		\end{equation}
		By \ref{F4}, there are $0<c_0\le C_0$ and $s_0>0$ such that $c_0|s|^{p-2}\le f(s)/s\le C_0|s|^{p-2}$ for $0<|s|\le s_0$. By integration, $F(s)\ge\frac{c_0}{p}|s|^p$ for $|s|\le s_0$. Using $\Psi(v)\le\lVert\n v\rVert_2^2$ on $\D$ and \eqref{scaling of trial function}, for $0<\t\le s_0/\lVert w_\mu\rVert_\infty$ we get
		\begin{equation*}
			\begin{aligned}
				I_\mu((w_\mu)_\t)&\le\t^N\Big(\lVert\n w_\mu\rVert_2^2+\frac{\mu\t^2}{2}\lVert w_\mu\rVert_2^2-\frac{c_0\t^p}{p}\irn w_\mu^p\,dx\Big)\\
				&=\mu^{\b+1}\t^N\Big(\lVert\n\omega\rVert_2^2+\frac{\t^2}{2}\lVert\omega\rVert_2^2-\frac{c_0\t^p}{p}\irn\omega^p\,dx\Big).
			\end{aligned}
		\end{equation*}
		Since $p>2$, there is $\t_1>0$, independent of $\mu$, such that the last bracket is negative for $\t\ge\t_1$. Moreover, by \eqref{scaling of trial function} we have $\t_1\le s_0/\lVert w_\mu\rVert_\infty$ for $\mu$ small, and therefore $I_\mu((w_\mu)_{\t_1})\le0$. {Since $I_\mu((w_\mu)_\t)>0$ for $0<\t\le\t_\mu(w_\mu)$, see Section \ref{sec action}, we get} $\t_\mu(w_\mu)<\t_1$ and
		\begin{equation*}
			a(\mu)\le\Lambda(\mu,w_\mu)=\max_{0<\t\le\t_1}I_\mu((w_\mu)_\t)\le\mu^{\b+1}\t_1^N\Big(\lVert\n\omega\rVert_2^2+\frac{\t_1^2}{2}\lVert\omega\rVert_2^2\Big)=C\mu^{\b+1}.
		\end{equation*}
		
		{\it Lower bound for $a(\mu)$ and the mass when $p\ge2+\frac4N$.} Let $u\in\cC(\mu)$. By \eqref{asymptotics at zero}, \eqref{energy identity on ground states} and \eqref{lowe bound of Psi}, 
		\begin{equation}\label{upper bound of L2gradient}
			C\mu^{\beta+1}\ge a(\mu)\ge\frac2{N+2}\Psi(u)\ge\frac1{N+2}\lVert\n u\rVert_2^2.
		\end{equation}
		Recall that $G(s):=F(s)+\frac1Nf(s)s$. The upper bound in \ref{F4} gives $G(s)\le C|s|^p$ for $|s|\le s_0$. By \eqref{upper bound of L2gradient} and \eqref{Linfty by L2}, we have $\lVert u\rVert_{\infty}\le s_0$ for $\mu$ small enough. Writing $\a:=\frac{N(p-2)}{4}$ and $\g:=\frac p2-\a$, the identity $M_\mu(u)=0$, \eqref{lowe bound of Psi} and the Gagliardo--Nirenberg inequality $\lVert u\rVert_p^p\le C\lVert\n u\rVert_2^{2\a}\lVert u\rVert_2^{2\g}$ give
		\begin{equation}\label{GN on ground states}
			\frac12 \lVert\n u\rVert_2^2 +\frac{N+2}{2N}\mu \lVert u\rVert_2^2\le\Psi(u)+\frac{N+2}{2N}\mu \lVert u\rVert_2^2=\irn G(u)\,dx\le C\irn|u|^p\,dx\le C\lVert\n u\rVert_2^{2\a}\lVert u\rVert_2^{2\g}.
		\end{equation}
		If $p\ge2+\frac4N$, then $\a\ge1$ and $\g>0$ since $p<2^*$. Therefore, \eqref{upper bound of L2gradient} and \eqref{GN on ground states} give 
		\begin{equation*}
			\lVert u\rVert_2^{2\g}\ge c\lVert\n u\rVert_2^{2(1-\a)}\ge c(C\mu^{\b+1})^{1-\a},
		\end{equation*}
		which gives the first inequality in \eqref{mass at zero}, because $(\b+1)(1-\a)/\g=\b$. In particular, when $p=2+4/N$, one has $\a=1$ and $\b=0$, and this estimate gives a uniform positive lower bound for the mass.
		
		Recall that $K(s):=\frac12f(s)s-F(s)$. We claim that $K(s)\ge c|s|^p$ for $|s|\le s_0$. For $\varepsilon\in\{-1,1\}$ let $q_\varepsilon(r):=\varepsilon f(\varepsilon r)/r$, $r>0$, which is increasing by \ref{F2}, and put {$\zeta:=(c_0/(2C_0))^{1/(p-2)}\in(0,1)$}. For $0<r\le s_0$ we have
		\[
		\begin{aligned}
			K(\varepsilon r)
			&=\int_0^r t\big(q_\varepsilon(r)-q_\varepsilon(t)\big)\,dt\\
			&\ge\int_0^{{\zeta} r}t\big(c_0r^{p-2}-C_0({\zeta} r)^{p-2}\big)\,dt
			=\frac{c_0{\zeta}^2}{4}r^p{,}
		\end{aligned}
		\]
		which proves the claim. Since $\lVert u\rVert_\infty\le s_0$, the energy identity \eqref{energy identity on ground states}, \eqref{GN on ground states} and the first inequality in \eqref{mass at zero} yield
		\begin{equation*}
			a(\mu)\ge\frac2{N+2}\irn K(u)\,dx\ge c\irn|u|^p\,dx\ge c\irn G(u)\,dx\ge c\mu\lVert u\rVert_2^2\ge c\mu^{\b+1},
		\end{equation*}
		which is the second inequality in \eqref{mass at zero}. In the critical case $p=2+4/N$, the upper and lower bounds therefore give $c\mu\le a(\mu)\le C\mu$ for all sufficiently small $\mu>0$.
		
		Since $\b+1>0$ if and only if $p<2^*$, and $\b<0$ if and only if $p>2+\frac4N$, assertions (i)--(iii) follow from \eqref{asymptotics at zero} and \eqref{mass at zero}.
	\end{proof}

	{Recall that $b_m(\mu)=a(\mu)-\mu m^2$, see \eqref{def of g_energy}. The next two lemmas collect the properties of $b_m$, $c_-$ and $a$ which are used in the proofs of the main results.}

	\begin{lemma}\label{properties of b}
		{Assume \ref{F1}--\ref{F4} and let $m>0$. Then $b_m$ is continuous on $(0,\infty)$,
		\begin{equation}\label{limits of b}
			\lim_{\mu\to0^+}b_m(\mu)=0
			\qquad\text{and}\qquad
			\lim_{\mu\to\infty}b_m(\mu)=+\infty.
		\end{equation}
		In particular, $\inf_{\mu>0}b_m(\mu)\le0$, and this infimum is attained whenever it is negative. Moreover, for every $\mu>0$ the following assertions hold.}
		\begin{enumerate}[label={\rm (\roman*)}, ref={\rm (\roman*)}]
			\item\label{b local min} {If $\mu$ is a local minimum point of $b_m$, then $c_-(\mu)=c_+(\mu)=m^2$, that is, $\cC(\mu)\subset\cS_m$.}
			\item\label{b normalized} {If $u\in\cC(\mu)\cap\cS_m$, then $(\mu,u)$ is a solution of \eqref{eq} with $\cE(u)=b_m(\mu)$. Moreover, $u$ has constant sign, it is radially symmetric about some point $x_0\in\RN$, and $u(\cdot+x_0)\in\cC(\mu)\cap\cS_m$. If $f$ is odd, then also $-u\in\cC(\mu)\cap\cS_m$.}
		\end{enumerate}
	\end{lemma}
	\begin{proof}
		{By Theorem \ref{relation btw eng' and mass}, $a$ has finite one-sided derivatives at every point, hence $a$ and $b_m$ are continuous on $(0,\infty)$. Since $b_m(\mu)=\mu\big(a(\mu)/\mu-m^2\big)$, the limits in \eqref{limits of b} follow from $a(\mu)\to0$ as $\mu\to0^+$ (Lemma \ref{behavior at zero}) and from $a(\mu)/\mu\to\infty$ as $\mu\to\infty$ (Lemma \ref{behavior at infinity}). Hence $\inf_{\mu>0}b_m(\mu)\le0$, and if this infimum is negative, then by \eqref{limits of b} it coincides with the infimum of $b_m$ over a compact subinterval of $(0,\infty)$, where it is attained.

		\ref{b local min} Since $\mu$ is a local minimum point of $b_m$, we have $b_m(\mu+h)\ge b_m(\mu)$ and $b_m(\mu-h)\ge b_m(\mu)$ for every sufficiently small $h>0$, that is,
		\begin{equation*}
			\frac{a(\mu+h)-a(\mu)}{h}\ge m^2\ge\frac{a(\mu)-a(\mu-h)}{h}.
		\end{equation*}
		Hence \eqref{upper bounds for one-sided derivatives}, with $\mu_0=\mu$, gives $c_-(\mu)\ge m^2\ge c_+(\mu)$. Since $c_-(\mu)\le c_+(\mu)$, we conclude that $c_-(\mu)=c_+(\mu)=m^2$.

		\ref{b normalized} Recall from Subsection \ref{subsec framework} that $u$ is a critical point of $I_\mu$. Since $u\in\cS_m$, the pair $(\mu,u)$ is a solution of \eqref{eq}, and
		\begin{equation*}
			\cE(u)=I_\mu(u)-\frac\mu2\lVert u\rVert_2^2=a(\mu)-\mu m^2=b_m(\mu).
		\end{equation*}
		The sign and the symmetry of $u$ are given by Theorem \ref{symmetry of action ground states}. Finally, if $f$ is odd, then $I_\mu$ and $M_\mu$ are even, hence $-u\in\cC(\mu)\cap\cS_m$.}
	\end{proof}

	\begin{lemma}\label{thresholds lemma}
		{Assume \ref{F1}--\ref{F4} with $2+\frac4N\le p<2^*$, and let
		\begin{equation*}
			c_*=\inf_{\mu>0}c_-(\mu),
			\qquad
			r_e=\inf_{\mu>0}\frac{a(\mu)}{\mu}
		\end{equation*}
		be as in Theorems \ref{main theorem A} and \ref{main theorem energy}. Then $c_*>0$, and the function $\mu\mapsto a(\mu)-c_*\mu$ is nonnegative and nondecreasing on $(0,\infty)$. In particular, $c_*\le r_e<\infty$. Moreover, if $p>2+\frac4N$, then $c_*$ and $r_e$ are attained, and $a(\mu)>c_*\mu$ for every $\mu>0$. In particular, $c_*<r_e$ in this case.}
	\end{lemma}
	\begin{proof}
		{Since $\b\le0$, by \eqref{mass at zero} there are $c,\delta>0$ such that $c_-(\mu)\ge c\mu^{\b}\ge c\delta^{\b}$ for $0<\mu\le\delta$. Moreover, \eqref{lower bounds at infinity} and the monotonicity of $\sigma$ give $c_-(\mu)\ge d_N\sigma(\delta)^{N+2}$ for $\mu\ge\delta$. Hence $c_*>0$. The function $g(\mu):=a(\mu)-c_*\mu$ is continuous on $(0,\infty)$ by Lemma \ref{properties of b}, and the second inequality in \eqref{upper bounds for one-sided derivatives}, together with $c_+(\mu)\ge c_-(\mu)\ge c_*$, gives
		\begin{equation*}
			\liminf_{h\to0^+}\frac{g(\mu)-g(\mu-h)}{h}\ge c_+(\mu)-c_*\ge0
			\qquad\text{for every }\mu>0.
		\end{equation*}
		Therefore $g$ is nondecreasing on $(0,\infty)$. Indeed, otherwise there are $\eps>0$ and $0<\mu_1<\mu_2$ such that $g_\eps(\mu):=g(\mu)+\eps\mu$ satisfies $g_\eps(\mu_1)>g_\eps(\mu_2)$. Then $\mu_3:=\min\{\mu\in[\mu_1,\mu_2]:g_\eps(\mu)\le g_\eps(\mu_2)\}>\mu_1$, and $g_\eps(\mu_3-h)>g_\eps(\mu_3)$ for $0<h\le\mu_3-\mu_1$, which is impossible, because the lower limit of the left difference quotients of $g_\eps$ at $\mu_3$ is not smaller than $\eps$. Since $g(\mu)\to0$ as $\mu\to0^+$ by Lemma \ref{behavior at zero}, $g$ is also nonnegative. Thus $a(\mu)/\mu\ge c_*$ for every $\mu>0$, and $c_*\le r_e<\infty$.

		Let now $p>2+\frac4N$. By Lemmas \ref{behavior at infinity} and \ref{behavior at zero}, $c_-(\mu)\to\infty$ and $a(\mu)/\mu\to\infty$ as $\mu\to0^+$ and as $\mu\to\infty$. Since $c_-$ is lower semicontinuous by Corollary \ref{attainment of c pm} and $\mu\mapsto a(\mu)/\mu$ is continuous, both infima are attained. Moreover, $a(\mu)-c_*\mu=\mu\big(a(\mu)/\mu-c_*\big)>0$ for $\mu$ small, and the monotonicity of $\mu\mapsto a(\mu)-c_*\mu$ gives $a(\mu)>c_*\mu$ for every $\mu>0$. In particular, $r_e=\min_{\mu>0}a(\mu)/\mu>c_*$.}
	\end{proof}
	
	\subsection{\texorpdfstring{{Proofs}}{Proofs} of Theorems \ref{main theorem A}--\ref{main theorem energy}}
	
	\begin{proof}[Proof of Theorem \ref{main theorem A}]
		{By Lemma \ref{thresholds lemma}, $c_*>0$, and $c_*$ is attained if $p>2+\frac4N$.}

		\ref{main A i} If $u\in\cC(\mu)$ has $\frac12\lVert u\rVert_2^2=m^2$, then $m^2\ge c_-(\mu)\ge c_*$, which excludes $m^2<c_*$. If $c_*$ is not attained, then $c_-(\mu)>c_*$ for every $\mu>0$, and the same inequality excludes $m^2=c_*$.

		\ref{main A ii} Let $m^2=c_*=c_-(\mu_0)$. {By Corollary \ref{attainment of c pm}, there is $u\in\cC(\mu_0)$ such that $\frac12\lVert u\rVert_2^2=c_-(\mu_0)=m^2$, that is, $u\in\cC(\mu_0)\cap\cS_m$. By Lemma \ref{properties of b}\ref{b normalized}, after a translation, $(\mu_0,u)$ is a radial solution of \eqref{eq} of constant sign, and
		\begin{equation*}
			\cE(u)=b_m(\mu_0)=a(\mu_0)-c_*\mu_0.
		\end{equation*}
		By Lemma \ref{thresholds lemma}, this quantity is nonnegative, and it is positive if $p>2+\frac4N$.}

		\ref{main A iii} Let $m^2>c_*$ and choose $\mu_1>0$ with $c_-(\mu_1)<m^2$. By Theorem \ref{relation btw eng' and mass}, $(b_m)'_+(\mu_1)=c_-(\mu_1)-m^2<0$, hence $b_m(\mu_1+\eta)<b_m(\mu_1)$ for some $\eta>0$. {Since $b_m$ is continuous and $b_m(\mu)\to+\infty$ as $\mu\to\infty$ by Lemma \ref{properties of b}, the minimum of $b_m$ over $[\mu_1,\infty)$ is attained at some $\mu>\mu_1$, which is a local minimum point of $b_m$. Lemma \ref{properties of b}\ref{b local min} then gives $c_-(\mu)=c_+(\mu)=m^2$, that is, $\frac12\lVert u\rVert_2^2=m^2$ for every $u\in\cC(\mu)$. By Lemma \ref{properties of b}\ref{b normalized}, after a translation, $(\mu,u)$ is a radial solution of \eqref{eq} of constant sign with $\cE(u)=b_m(\mu)$.

		Finally, if $f$ is odd, then in \ref{main A ii} and \ref{main A iii} we may replace $u$ by $-u$, see Lemma \ref{properties of b}\ref{b normalized}, so that the solutions can be taken positive. }
	\end{proof}
	
	\begin{remark}
		{Let $2+\frac4N<p<2^*$ and $m^2>c_*$, and let $\mu_1<\mu$ be as in the proof of Theorem \ref{main theorem A}\ref{main A iii}. Then $b_m$ also has a local maximum point $\bar\mu\in(0,\mu)$. Indeed, $b_m(\nu)=\nu\big(a(\nu)/\nu-m^2\big)>0$ for $\nu$ small by Lemma \ref{behavior at zero}, $b_m(\nu)\to0$ as $\nu\to0^+$ by \eqref{limits of b}, and $b_m(\mu)<b_m(\mu_1)$. Hence the maximum of $b_m$ over $(0,\mu]$ is positive, and it is attained at some $\bar\mu\in(0,\mu)$.} In contrast with a local minimum point, Theorem \ref{relation btw eng' and mass} only gives $c_-(\bar\mu)\le m^2\le c_+(\bar\mu)$ at $\bar\mu$, which does not guarantee that some element of $\cC(\bar\mu)$ has the prescribed mass. Whether $\bar\mu$ gives rise to a normalized solution {is an open question}.
	\end{remark}
	
	\begin{proof}[Proof of Theorem \ref{main theorem sub}]
		{Since $a(\mu)/\mu\to0$ as $\mu\to0^+$ by Lemma \ref{behavior at zero}, we have $b_m(\mu)=\mu\big(a(\mu)/\mu-m^2\big)<0$ for $\mu$ small. Hence, by Lemma \ref{properties of b}, $b_m$ attains a negative global minimum at some $\mu>0$, and Lemma \ref{properties of b}\ref{b local min} gives $c_-(\mu)=c_+(\mu)=m^2$, that is, $\frac12\lVert u\rVert_2^2=m^2$ for every $u\in\cC(\mu)$. By Lemma \ref{properties of b}\ref{b normalized}, after a translation, $(\mu,u)$ is a radial solution of \eqref{eq} of constant sign with
		\begin{equation*}
			\cE(u)=b_m(\mu)=\min_{(0,\infty)}b_m<0,
		\end{equation*}
		and if $f$ is odd, we may replace $u$ by $-u$, so that the solution can be taken positive.}
	\end{proof}
	
	\begin{proof}[Proof of Theorem \ref{main theorem energy}]
		{{\it Step 1: $e_m\le b_m(\mu)$ for every $\mu>0$.} Let $\mu>0$, $u\in\cC(\mu)$ and $\t:=\big(2m^2/\lVert u\rVert_2^2\big)^{\frac1{N+2}}$, so that $u_\t\in\cS_m$. Since $u\in\cM_\mu$, \eqref{strict maximum on fibers} gives
		\begin{equation}\label{dilation onto the sphere}
			e_m\le\cE(u_\t)=I_\mu(u_\t)-\mu m^2\le I_\mu(u)-\mu m^2=b_m(\mu),
		\end{equation}
		and the second inequality is strict unless $\t=1$, that is, unless $u\in\cS_m$.}

		{\it Step 2: $\cE(u)\ge b_m(\mu_u)$ for every $u\in\cS_m$ with $\cE(u)\le0$, where}
		\begin{equation*}
			\mu_u:=\frac{N}{(N+2)m^2}\Big(\irn G(u)\,dx-\Psi(u)\Big),\qquad G(s)=F(s)+\frac1Nf(s)s ,
		\end{equation*}
		{{\it and equality holds if and only if $u\in\cC(\mu_u)$.}} Indeed, $\cE(u)\le0$ gives $\irn F(u)\,dx\ge\Psi(u)$, and $\irn f(u)u\,dx>0$ because $f(s)s>0$ for $s\ne0$ by \ref{F1}--\ref{F2} and $u\ne0$. Thus $\mu_u>0$. Since $\lVert u\rVert_2^2=2m^2$, by \eqref{function Hmu}
		\begin{equation*}
			M_{\mu_u}(u)=\Psi(u)-\irn G(u)\,dx+\mu_u\Big(\frac12+\frac1N\Big)\lVert u\rVert_2^2=\Psi(u)-\irn G(u)\,dx+\frac{(N+2)m^2}{N}\,\mu_u=0,
		\end{equation*}
		i.e. $u\in\cM_{\mu_u}$. Hence $I_{\mu_u}(u)\ge a(\mu_u)$ and $\cE(u)=I_{\mu_u}(u)-\mu_um^2\ge a(\mu_u)-\mu_um^2=b_m(\mu_u)${, with equality if and only if $I_{\mu_u}(u)=a(\mu_u)$, that is, $u\in\cC(\mu_u)$}.

		{{\it Step 3: proof of \ref{energy i} and \ref{energy ii}.} By Step 1 and since $b_m(\mu)\to0$ as $\mu\to0^+$ by Lemma \ref{behavior at zero}, we have $e_m\le\inf_{\mu>0}b_m(\mu)\le0$. Conversely, let $u\in\cS_m$. If $\cE(u)>0$, then $\cE(u)>\inf_{\mu>0}b_m(\mu)$, and if $\cE(u)\le0$, then $\cE(u)\ge b_m(\mu_u)\ge\inf_{\mu>0}b_m(\mu)$ by Step 2. Hence
		\begin{equation}\label{Legendre identity}
			e_m=\inf_{\mu>0}b_m(\mu)\le0,
		\end{equation}
		where $e_m>-\infty$ by \eqref{lower bound of energy}, and $m^2\mapsto e_m$ is concave and nonincreasing as the infimum of the affine nonincreasing functions $m^2\mapsto a(\mu)-\mu m^2$. This proves \ref{energy i}.

		Let $\mu_m$ be a global minimum point of $b_m$ and let $u\in\cC(\mu_m)$. Then $b_m(\mu_m)=e_m$ by \eqref{Legendre identity}, so that equality holds throughout \eqref{dilation onto the sphere} with $\mu=\mu_m$. Hence $u\in\cS_m$ and $\cE(u)=e_m$. Thus $\cC(\mu_m)\subset\cS_m$, that is, $c_-(\mu_m)=c_+(\mu_m)=m^2$, and every element of $\cC(\mu_m)$ is an energy ground state. Conversely, let $u\in\cS_m$ be such that $\cE(u)=e_m$. Since $e_m\le0$, Step 2 and \eqref{Legendre identity} give
		\begin{equation*}
			e_m=\cE(u)\ge b_m(\mu_u)\ge\inf_{\mu>0}b_m(\mu)=e_m.
		\end{equation*}
		Hence $\mu_u$ is a global minimum point of $b_m$, and $u\in\cC(\mu_u)$ by the equality case in Step 2. In particular, $e_m$ is attained if and only if $b_m$ attains its infimum on $(0,\infty)$. The last assertion of \ref{energy ii} follows from Lemma \ref{properties of b}\ref{b normalized}.

		{\it Step 4: proof of \ref{energy iii}.} If $2<p<2+\frac4N$, then $b_m$ attains a negative global minimum, as shown in the proof of Theorem \ref{main theorem sub}, and \ref{energy i}--\ref{energy ii} give the claim. Let $2+\frac4N\le p<2^*$. By Lemma \ref{thresholds lemma}, $c_*\le r_e<\infty$, and $c_*<r_e$ if $p>2+\frac4N$. Recall that $b_m(\mu)=\mu\big(a(\mu)/\mu-m^2\big)$.
		If $m^2<r_e$, then $b_m>0$ on $(0,\infty)$. Hence, by \eqref{limits of b}, $\inf_{\mu>0}b_m(\mu)=0$ is not attained, and $e_m=0$ is not attained by \ref{energy i}--\ref{energy ii}.
		If $m^2>r_e$, then $b_m$ is negative somewhere. Hence $\inf_{\mu>0}b_m(\mu)<0$ is attained by Lemma \ref{properties of b}, and $e_m<0$ is attained.
		If $m^2=r_e$, then $b_m\ge0$ on $(0,\infty)$, so that $e_m=0$, and $b_m(\mu)=0$ if and only if $a(\mu)/\mu=r_e$. By \ref{energy ii}, $e_m=0$ is attained if and only if $\mu\mapsto a(\mu)/\mu$ attains its infimum, which is the case if $p>2+\frac4N$ by Lemma \ref{thresholds lemma}.}
	\end{proof}
	
	\begin{remark}
		For $m^2>r_e$, and for $m^2=r_e$ when the infimum $r_e$ is attained, Theorem \ref{main theorem energy} provides a global minimum point $\mu_m$ of $b_m$, and every action ground state at the frequency $\mu_m$ is an energy ground state of mass $2m^2$.
		The frequency $\mu$ selected in Theorem \ref{main theorem A}\ref{main A iii}, however, minimizes $b_m$ only on a half-line $[\mu_1,\infty)$ with $\mu_1<\mu$, as in its proof, and the above arguments do not tell whether this minimum is a global one on $(0,\infty)$. Thus the elements of $\cC(\mu)$ are energy ground states if and only if $b_m(\mu)=\inf_{\nu>0}b_m(\nu)$.
		For $c_*<m^2<r_e$ they have positive energy, so they are not energy ground states, and in fact $e_m=0$ is not attained.
		At $m^2=c_*$, when $c_*$ is attained, Theorem \ref{main theorem A}\ref{main A ii} provides the elements of $\cC(\mu_0)$ of mass $2c_-(\mu_0)$, but it does not assert that every element of $\cC(\mu_0)$ has this mass.
		The correspondence between the energy ground states and the global minimum points of $b_m$ is the counterpart of the rigidity result \cite[Theorem 1.3]{DoSeTi}. {Finally, we point out that Steps 1--3 of the proof of Theorem \ref{main theorem energy}, and hence assertions \ref{energy i} and \ref{energy ii} of that theorem, use neither Lemma \ref{compactness lemma} nor Theorem \ref{relation btw eng' and mass}.}
	\end{remark}
	
	\begin{proof}[Proof of Corollary \ref{critical power corollary}]
		Put $p_c:=2+4/N$ and define
		\[
		\cE_\Delta(v):=\frac12\lVert\nabla v\rVert_2^2
		-\frac1{p_c}\lVert v\rVert_{p_c}^{p_c},
		\qquad v\in H^1(\RN).
		\]
		The sharp Gagliardo--Nirenberg inequality \cite{Weinstein} states that
		\[
		\lVert v\rVert_{p_c}^{p_c}
		\le\frac{p_c}{2\lVert Q\rVert_2^{4/N}}
		\lVert\nabla v\rVert_2^2\lVert v\rVert_2^{4/N}
		\qquad\text{for every }v\in H^1(\RN).
		\]
		Since $1-\sqrt{1-r^2}>r^2/2$ for $0<r\le1$ and a nonzero
		$H^1$ function cannot have zero gradient almost everywhere,
		for every $u\in\cS_m$ we obtain
		\[
		\cE(u)>\cE_\Delta(u)
		\ge\frac12\left[1-\left(\frac{m}{m_c}\right)^{4/N}\right]
		\lVert\nabla u\rVert_2^2.
		\]
		Thus $\cE(u)>0$ for $0<m\le m_c$. Since $e_m\le0$ by {\eqref{Legendre identity}},
		it follows that
		$e_m=0$ and that the infimum is not attained in this range.
		
		Let now $m>m_c$ and put $\alpha:=m/m_c>1$. The Nehari and
		Poho\v{z}aev identities for $Q$ give
		$\lVert Q\rVert_{p_c}^{p_c}=(p_c/2)\lVert\nabla Q\rVert_2^2$.
		Consequently, $\lVert\alpha Q\rVert_2^2=2m^2$ and
		\[
		\cE_\Delta(\alpha Q)
		=\frac{\alpha^2-\alpha^{p_c}}2\lVert\nabla Q\rVert_2^2<0.
		\]
		Approximating $\alpha Q$ in $H^1$ by smooth compactly supported
		functions and normalizing their $L^2$-norms, we may therefore
		choose $v\in\cC_c^\infty(\RN)$ with
		$\lVert v\rVert_2^2=2m^2$ and $\cE_\Delta(v)<0$.
		For $\lambda>0$ set $v_\lambda(x)=\lambda^{N/2}v(\lambda x)$.
		For all sufficiently small $\lambda$,
		\[
		\lVert v_\lambda\rVert_2^2=2m^2,\qquad
		\lVert\nabla v_\lambda\rVert_\infty
		=\lambda^{N/2+1}\lVert\nabla v\rVert_\infty<1,
		\]
		hence $v_\lambda\in\cS_m$. Rationalizing the kinetic density
		and changing variables gives
		\[
		\frac{\Psi(v_\lambda)}{\lambda^2}
		=\int_{\RN}\frac{|\nabla v|^2}
		{1+\sqrt{1-\lambda^{N+2}|\nabla v|^2}}\,dx
		\longrightarrow\frac12\lVert\nabla v\rVert_2^2
		\quad\text{as }\lambda\to0^+.
		\]
		The integrand is bounded by $|\nabla v|^2$, so the limit follows
		from dominated convergence. Since
		$\lVert v_\lambda\rVert_{p_c}^{p_c}
		=\lambda^2\lVert v\rVert_{p_c}^{p_c}$, we obtain
		\[
		\lambda^{-2}\cE(v_\lambda)\longrightarrow\cE_\Delta(v)<0.
		\]
		Hence $e_m<0$ for every $m>m_c$. The sign characterization in
		Theorem \ref{main theorem energy}\ref{energy iii} now gives
		$r_e=m_c^2$ and attainment for $m>m_c$. Since $e_{m_c}=0$
		is not attained, the same theorem shows that the infimum of
		$a(\mu)/\mu$ is not attained either.
		
		As for the comparison with the Laplacian, the same sharp inequality and
		the mass-preserving dilations $v_\lambda$ show that the infimum of $\cE_\Delta$ on the mass sphere in $H^1(\RN)$ is zero for $0<m\le m_c$. It is not attained for $m<m_c$, by the displayed lower bound for $\cE_\Delta$,
		and it is attained by $Q$ at $m=m_c$. For $m>m_c$, a function $v$ with negative energy
		as above satisfies
		$\cE_\Delta(v_\lambda)=\lambda^2\cE_\Delta(v)\to-\infty$
		as $\lambda\to\infty$. These concentrating dilations are admissible
		on the mass sphere in $H^1(\RN)$, but not in $\cS_m$.
	\end{proof}
	
	\appendix
	\section{Two estimates for functions in $\D$}\label{sec estimates}
	
	In this appendix we prove the two estimates announced in Section \ref{sec intro}, which express the fact that the gradient constraint $\lVert\n u\rVert_\infty\le1$ prevents concentration.
	
	\begin{lemma}\label{Linfty bound}
		Let $u\in\D$. Then
		\begin{equation}\label{Linfty by L2}
			\lVert u\rVert_\infty\le2\left(\frac{\lVert u\rVert_2^2}{|B_1|}\right)^{\frac1{N+2}},
			\qquad
			\lVert u\rVert_\infty\le C_N'\lVert\n u\rVert_2^{2/N},
		\end{equation}
		where $C_N'$ depends only on $N$.
	\end{lemma}
	\begin{proof}
		By \eqref{embedding}, $u\in\cC_0(\RN)$, hence $A:=\lVert u\rVert_\infty=|u(x_0)|$ for some $x_0\in\RN$, and we may assume $A>0$. Since $u$ is $1$-Lipschitz, $|u|\ge A/2$ on $B_{A/2}(x_0)$, hence
		\begin{equation*}
			\lVert u\rVert_2^2\ge\int_{B_{A/2}(x_0)}|u|^2\,dx\ge|B_1|\Big(\frac A2\Big)^{N+2},
		\end{equation*}
		which is the first inequality in \eqref{Linfty by L2}. For the second one, set $w=(|u|-A/2)_+$. Then $w\in H^1(\RN)$, $\lVert\n w\rVert_2\le\lVert\n u\rVert_2$ and $w\ge A/4$ on $B_{A/4}(x_0)$. The Sobolev inequality then gives
		\begin{equation*}
			\frac A4|B_1|^{1/2^*}\Big(\frac A4\Big)^{N/2^*}\le\lVert w\rVert_{2^*}\le S_N\lVert\n u\rVert_2 .
		\end{equation*}
		Since $1+N/2^*=N/2$, we get $A^{N/2}\le C\lVert\n u\rVert_2$, and the second inequality in \eqref{Linfty by L2} follows.
	\end{proof}
	
	\begin{lemma}\label{energy bounded below}
		Assume \ref{F1}--\ref{F2} and let $m>0$. Then
		\begin{equation}\label{lower bound of energy}
			\inf_{u\in\cS_m}\cE(u)\ge-\kappa(S(m))\,m^2>-\infty,
		\end{equation}
		where $S(m):=2\big(2m^2/|B_1|\big)^{\frac1{N+2}}$ and $\kappa(S):=\max\big\{f(S)/S,\,f(-S)/(-S)\big\}$.
	\end{lemma}
	\begin{proof}
		Let $u\in\cS_m$. By \eqref{Linfty by L2}, $\lVert u\rVert_\infty\le S(m)$. Since $s\mapsto f(s)/s$ is positive, nondecreasing on $(0,\infty)$ and nonincreasing on $(-\infty,0)$ by \ref{F1}--\ref{F2}, we have $0<f(s)/s\le\kappa(S)$ for $0<|s|\le S$, hence
		\begin{equation}\label{quadratic bound for F}
			F(s)=\int_0^s\frac{f(t)}{t}\,t\,dt\le\kappa(S)\cdot\frac{s^2}{2}\qquad\text{for every }S>0\text{ and }|s|\le S .
		\end{equation}
		Applied with $S=S(m)$, this gives $\irn F(u)\,dx\le\frac12\kappa(S(m))\lVert u\rVert_2^2=\kappa(S(m))m^2$, and \eqref{lower bound of energy} follows from $\Psi\ge0$.
	\end{proof}
	
	\section{Symmetry and monotonicity of the action ground states}\label{sec ground states}
	
	In this appendix we prove that every action ground state of \eqref{eq-mu} has constant sign and is radially symmetric about some point, with a strictly decreasing radial profile of its absolute value, see Theorem \ref{symmetry of action ground states}. In the application to the normalized problem, the symmetry of {\it every} action ground state is needed, and not only the existence of a radial minimizer at each frequency: for instance, at $m^2=c_*$ we first select a ground state of minimum mass and then center it by a translation, which does not change its mass, and the same argument gives the radial symmetry of every energy ground state. 
	
	We combine the Schwarz rearrangement with the equality case of the P\'olya--Szeg\H{o} inequality due to Brothers and Ziemer \cite{BZ}, in the form given by Burchard and Ferone \cite{BF}, and with the uniqueness for the radial ordinary differential equation. The point is that the Born--Infeld density is a series of even powers of $|\nabla u|$ with positive coefficients, and hence the equality $\Psi(|u|^*)=\Psi(u)$, which follows from the minimality, reduces to the equality $\lVert\nabla|u|^*\rVert_2=\lVert\nabla u\rVert_2$. Throughout, $N\ge3$, $\mu>0$ and $f$ satisfies \ref{F1}--\ref{F3}, and we use the notation of Subsection \ref{subsec framework}. {For a nonnegative $v\in\D$} we denote by $v^*$ its Schwarz rearrangement, which is nonnegative, radial, radially nonincreasing and equimeasurable with $v$, see \cite[Chapter 3]{LL}.
	
	We recall from Subsection \ref{subsec framework} that every $u\in\cC(\mu)$ is a critical point of $I_\mu$ in the sense of \eqref{def critical point}, hence a solution of \eqref{eq-mu} of class $\cC^2(\RN)$ with $\lVert\nabla u\rVert_\infty<1$, see \cite[Lemmas 2.3 and 5.5]{BIM}.
	
	\begin{lemma}\label{constant sign lemma}
		Every $u\in\cC(\mu)$ has constant sign, i.e., $u>0$ or $u<0$ a.e. on $\RN$.
	\end{lemma}
	\begin{proof}
		See \cite[Proof of Theorem 1.6]{BIM}.
	\end{proof}
	
	\begin{lemma}\label{rearrangement lemma}
		Let $u\in\cC(\mu)$ be nonnegative. Then $u^*\in\cC(\mu)$, $\Psi(u^*)=\Psi(u)$ and $\lVert\nabla u^*\rVert_2=\lVert\nabla u\rVert_2$.
	\end{lemma}
	\begin{proof}
		{By \cite[Proof of Theorem 1.6]{BIM}, applied to the equation with unit frequency and with $f$ replaced by $g_\mu$ as in Subsection \ref{subsec framework}, we have $u^*\in\cC(\mu)$. As in that proof, see \cite[Proof of Theorem 1.2]{BIM}, the P\'olya--Szeg\H{o} inequality applies to every term of the expansion $1-\sqrt{1-r^2}=\sum_{k=1}^\infty a_kr^{2k}$, $0\le r\le1$, where $a_k>0$ and $a_1=\frac12$. Since $|\nabla u|\le1$ and $|\nabla u^*|\le1$ a.e. on $\RN$, the monotone convergence theorem gives
		\begin{equation}\label{Delta series}
			\Psi(u)-\Psi(u^*)
			=\sum_{k=1}^\infty a_k\Big(\irn|\nabla u|^{2k}\,dx-\irn|\nabla u^*|^{2k}\,dx\Big),
		\end{equation}
		where every term of the series is nonnegative, see also \cite[Lemma 9]{ByeonJeanjeanMaris}. On the other hand, since $M_\mu(u)=M_\mu(u^*)=0$ and since $u$ and $u^*$ are equimeasurable,
		\begin{equation*}
			\Psi(u^*)=\irn H_\mu(u^*)\,dx=\irn H_\mu(u)\,dx=\Psi(u).
		\end{equation*}
		Hence all the terms in \eqref{Delta series} vanish, and the term with $k=1$ gives $\lVert\nabla u^*\rVert_2=\lVert\nabla u\rVert_2$.}
	\end{proof}
	
	\begin{lemma}\label{radial profile lemma}
		Let $v\in\cC(\mu)$ and $\varepsilon\in\{-1,1\}$ be such that $\varepsilon v(x)=W(|x|)$ for $x\in\RN$, where $W\ge0$ is nonincreasing on $[0,\infty)$. Then $W\in\cC^2([0,\infty))$, $W(r)>0$ for every $r\ge0$ and $W'(r)<0$ for every $r>0$.
	\end{lemma}
	\begin{proof}
		Since $v\in\cC^2(\RN)$ is a solution of \eqref{eq-mu} with $\lVert\nabla v\rVert_\infty<1$, the function $W(r)=\varepsilon v(re_1)$ belongs to $\cC^2([0,\infty))$, $W'(0)=0$, $|W'|\le\lVert\nabla v\rVert_\infty<1$ and $W'\le0$ on $[0,\infty)$. Put $h_\varepsilon(s):=\varepsilon f(\varepsilon s)-\mu s$, which is locally Lipschitz continuous by \ref{F1} with $h_\varepsilon(0)=0$. Writing the equation for $\varepsilon v$ in polar coordinates, for $r>0$ the profile $W$ solves
		\begin{equation}\label{radial ODE}
			W''=-\frac{N-1}{r}\,W'\,(1-W'^2)-h_\varepsilon(W)\,(1-W'^2)^{3/2}.
		\end{equation}
		Since the right-hand side of \eqref{radial ODE} is locally Lipschitz continuous in $(W,W')$ for $|W'|<1$ and $r>0$, the solutions of \eqref{radial ODE} with $|W'|<1$ are uniquely determined by their data at any $r_0>0$.
		
		If $W(r_0)=0$ for some $r_0>0$, then $W=0$ on $[r_0,\infty)$ because $W$ is nonnegative and nonincreasing, and $W\equiv0$ is the solution of \eqref{radial ODE} with the data $(0,0)$ at $r_0$, because $h_\varepsilon(0)=0$. By uniqueness, $W=0$ on $(0,\infty)$, which contradicts $v\neq0$. Hence $W>0$ on $(0,\infty)$, and $W(0)\ge W(r)>0$ by monotonicity.
		
		If $W'(r_0)=0$ for some $r_0>0$, then $r_0$ is a maximum point of $W'\le0$, hence $W''(r_0)=0$, and \eqref{radial ODE} gives $h_\varepsilon(W(r_0))=0$. The constant function $W(r_0)$ then solves \eqref{radial ODE} with the same data at $r_0$, and uniqueness forces $W$ to be constant on $(0,\infty)$, which contradicts the positivity of $W$ and $W(r)\to0$ as $r\to\infty$. Hence $W'<0$ on $(0,\infty)$.
	\end{proof}
	
	We next recall the equality case of the P\'olya--Szeg\H{o} inequality, due to Brothers and Ziemer \cite[Theorem 1.1]{BZ}, in the quantitative form of Burchard and Ferone \cite{BF}. {As in \cite{BF}, a nonnegative measurable function $z$ on $\RN$ is said to vanish at infinity if the level sets $\{z>t\}$, $t>0$, have finite measure, and $\omega_k$ denotes the Lebesgue measure of the unit ball in $\R^k$.}
	
	\begin{lemma}[{\cite[Theorem 1.1]{BF}}]\label{BZ lemma}
		{Let $1<p<\infty$, and let $z$ be a nonnegative function on $\RN$ which vanishes at infinity and whose distributional gradient belongs to $L^p(\RN)$. If $\lVert\nabla z\rVert_p=\lVert\nabla z^*\rVert_p$, where $z^*$ is the Schwarz rearrangement of $z$, then there exists a translation $\tau$ of $\RN$ such that
		\begin{equation}\label{BF estimate}
			\lVert z-z^*\circ\tau\rVert_q\le K_N^{1/q}\,\lVert z\rVert_q^{1/N'}\,\lVert z\chi_C\rVert_q^{1/N}
		\end{equation}
		for every $q\ge1$ with $z\in L^q(\RN)$. Here $N'=N/(N-1)$, $K_N=2\omega_{N-1}/\omega_N$ and}
		\begin{equation*}
			{C=\big\{x\in\RN:0<z(x)<\operatorname{ess\,sup}z\ \text{ and }\ |\nabla z(x)|=0\big\}.}
		\end{equation*}
	\end{lemma}

	{As observed in \cite{BF} right after the statement of this theorem, see also \cite[Lemma 3.3]{BF}, under the assumptions of Lemma \ref{BZ lemma} the function $z$ may be replaced by $z^*$ on the right-hand side of \eqref{BF estimate} and in the definition of $C$. In particular, if $z\in L^q(\RN)$ for some $q\ge1$ and the set
	\begin{equation*}
		C^*:=\big\{x\in\RN:0<z^*(x)<\operatorname{ess\,sup}z\ \text{ and }\ |\nabla z^*(x)|=0\big\}
	\end{equation*}
	has measure zero, then $z=z^*\circ\tau$ a.e. on $\RN$.}
	
	\begin{theorem}\label{symmetry of action ground states}
		Let $N\ge3$, assume \ref{F1}--\ref{F3}, and let $\mu>0$. Then every $u\in\cC(\mu)$ has the form
		\[
		u=\varepsilon\,W(|\cdot-x_0|),\qquad\varepsilon\in\{-1,1\},\quad x_0\in\RN,
		\]
		where $W\in\cC^2([0,\infty))$, $W(r)>0$ for $r\ge0$ and $W'(r)<0$ for $r>0$. Moreover, $u(\cdot+x_0)\in\cC(\mu)\cap X_r$ and $\lVert u(\cdot+x_0)\rVert_2=\lVert u\rVert_2$.
	\end{theorem}
	\begin{proof}
		Let $u\in\cC(\mu)$. By Lemma \ref{constant sign lemma}, we may assume that $u>0$ and prove the theorem with $\varepsilon=1$ (if $u<0$, the same argument applies to $-u$ and to the nonlinearity $s\mapsto-f(-s)$). By Lemma \ref{rearrangement lemma}, $u^*\in\cC(\mu)$ and
		\begin{equation}\label{equality of Dirichlet norms}
			\lVert\nabla u^*\rVert_2=\lVert\nabla u\rVert_2.
		\end{equation}
		Write $u^*(x)=W(|x|)$, where $W\ge0$ is nonincreasing. Lemma \ref{radial profile lemma} applies to $u^*$ with $\varepsilon=1$ and gives
		\begin{equation}\label{profile properties}
			W\in\cC^2([0,\infty)),\qquad W(r)>0\ \hbox{ for } r\ge0,\qquad W'(r)<0\ \hbox{ for } r>0.
		\end{equation}
		
		We now show that $u$ is a translate of $u^*$ by means of Lemma \ref{BZ lemma}, applied to $z=u$ {with $p=q=2$}. Its hypotheses hold: $u\in L^2(\RN)$ is nonnegative{, hence it vanishes at infinity}, $\nabla u\in L^2(\RN)$, and $\lVert\nabla u\rVert_2=\lVert\nabla u^*\rVert_2$ by \eqref{equality of Dirichlet norms}. Moreover, $\operatorname{ess\,sup}u=\lVert u^*\rVert_\infty=W(0)$, and since $W$ is strictly decreasing by \eqref{profile properties},
		\begin{equation*}
			\{x\in\RN:0<u^*(x)<W(0)\}=\RN\setminus\{0\}
		\end{equation*}
		and $\nabla u^*(x)=W'(|x|)\,x/|x|\neq0$ for $x\neq0$. Hence {the set $C^*$ associated with $z=u$ is empty, and the observation following Lemma \ref{BZ lemma} gives a translation $\tau$ such that $u=u^*\circ\tau$ a.e. on $\RN$. Writing $\tau(x)=x-x_0$ with $x_0\in\RN$, we get} $u=u^*(\cdot-x_0)$ a.e. on $\RN$. Both sides are continuous, hence $u(x)=u^*(x-x_0)=W(|x-x_0|)>0$ for every $x\in\RN$, with $W$ as in \eqref{profile properties}.
		
		Finally, the action, the constraint and the {$L^2$-norm} are invariant under translations, hence $u(\cdot+x_0)=W(|\cdot|)\in\cC(\mu)\cap X_r$ and $\lVert u(\cdot+x_0)\rVert_2=\lVert u\rVert_2$.
	\end{proof}
	
	The normalized action ground states of Theorems \ref{main theorem A} and \ref{main theorem sub} and the energy ground states of Theorem \ref{main theorem energy} therefore have the profile described in Theorem \ref{symmetry of action ground states}, and centering them does not change their mass. More generally,
	\[
	\big\{\lVert u\rVert_2^2:u\in\cC(\mu)\big\}=\big\{\lVert u\rVert_2^2:u\in\cC(\mu)\cap X_r\big\}\qquad\hbox{for every }\mu>0,
	\]
	which justifies the use of centered radial representatives in the analysis of the masses.
	
		\begin{remark}
		Assumption \ref{F4} is not used in this appendix. {The symmetry and strict monotonicity argument also applies to the zero mass problem of \cite{BIM}, under assumptions (F1)--(F6) of that paper, with $\mu=0$ and the space $X_0\subset L^{2^*}(\RN)$ of \cite{BIM} in place of $X$: it suffices to replace the results of \cite{BIM} quoted above by \cite[Theorem 1.2 and its proof, Lemmas 2.3 and 4.6]{BIM}, and to apply Lemma \ref{BZ lemma} with $p=2$ and $q=2^*$.}
	\end{remark}

	\section*{Acknowledgements}
	J. M. was partly supported by the National Science Centre, Poland, under the research grant 2025/57/B/ST1/04804 and the research grant 2023/51/B/ST1/00968.
	
	\smallskip
	\noindent\textbf{Data availability.} The manuscript has no associated data.
	
	\smallskip
	\noindent\textbf{Conflict of interest.} The authors declare no conflict of interest.

\end{document}